%% file: Wolpear.tex
\documentclass[a4paper]{amsart}
\usepackage[dvipsnames]{xcolor}
\usepackage{amsmath,amssymb,amsthm,mathrsfs,mathtools}
\usepackage{graphicx}
\usepackage{enumerate,enumitem} 
\usepackage{hyperref}
 \hypersetup{
     citecolor=CadetBlue,
     colorlinks=true,
     linkcolor=Periwinkle,
     filecolor=Periwinkle,      
     urlcolor=black,
     }
\usepackage{fancyhdr}
\usepackage{overpic}
\usepackage{orcidlink}
\usepackage{faktor}
\usepackage{xfrac}

\DeclareRobustCommand{\SkipTocEntry}[5]{}
\input{commands.tex}

\title{A Wolpert Theorem for shearing coordinates}
\author{Marie Abadie \orcidlink{0000-0002-8951-4643}}
\email{marie.abadie@uni.lu}
\date{\today}

\begin{document}

\begin{abstract}
   A theorem from Wolpert bounds the Weil--Petersson distance from a hyperbolic surface to a maximally noded surface in terms of the lengths of the underlying pants decomposition. We prove analogues for punctured surfaces both in terms of the lambda-lengths and in terms of the shearing coordinates of an ideal triangulation.
\end{abstract}
  
\maketitle
\tableofcontents

\section{Introduction}

The Weil--Petersson (WP) metric is a negatively curved Kähler metric on the Teichmüller space $\teich_{g,n}$ of a surface of genus $g$ with $n$ punctures~\cite{Ahlfors1961,Tromba1986,Wolpert1986a}. It is not complete~\cite{Wolpert1975,Chu1976}, and its completion is the augmented Teichmüller space $\overline{\teich}_{g,n}$, obtained by adding strata of noded surfaces~\cite{Bers1974,Abikoff1977,Masur1976}. The Weil--Petersson distance between two points of Teichmüller space is difficult to compute. Nevertheless, Wolpert~\cite{Wolpert2003,Wolpert2008} gave a simple upper bound for the distance from a point $X$ to the stratum $S(\Gamma)$ obtained by pinching a multicurve $\Gamma$:
\[
    \dwp(X,S(\Gamma)) \leq \sqrt{2\pi \sum_{\gamma \in \Gamma} \ell_X(\gamma)}.
\]
When $\Gamma$ is a pants decomposition, the stratum $S(\Gamma)$ is a single maximally noded surface, and the bound depends only on the Fenchel--Nielsen lengths with respect to $\Gamma$. Combined with the existence of short pants decompositions~\cite{Bers1974, Parlier2023}, it shows that every point of $\overline{\teich}_{g,n}$ lies at distance bounded in terms of $g$ and $n$ from a maximally noded surface. This result is used to study the coarse geometry of the WP metric, see for instance~\cite{Brock2003,CavendishParlier2012}. In particular Cavendish and Parlier obtained an upper bound of order $\sqrt{g}\log g$ for the WP diameter of moduli space with a fixed number of punctures~\cite{CavendishParlier2012}, which was recently shown to be sharp by Wu~\cite{wu2026}.

Punctured surfaces admit another natural decomposition, into ideal triangles. Whereas a pair of pants is determined up to isometry by its three boundary lengths, all ideal triangles are isometric, so the geometry of the surface is entirely captured by the way the triangles are glued together. Their gluing is determined by the shearing coordinates $s_e(X)$ associated to each arc $e$ of the triangulation~\cite{Bonahon1996,fock1998,Thurston1998}. Bestvina, Bromberg, Fujiwara and Souto showed that the length of every essential closed curve is a convex function of the shearing coordinates~\cite[Theorem~1.3]{BestvinaBrombergFujiwaraSouto2013}. We want an explicit bound for the WP distance between two surfaces in terms of their shearing coordinates.

We first establish a distance estimate in terms of truncated arc-lengths. After truncating the cusps along horocycles of length $1/4$, each arc $e$ of an ideal triangulation $\T$ of a hyperbolic surface $X$ has a finite \emph{truncated length} $\delta_e(X)$, see Section~\ref{subsec:truncation}. We write
\[
    \mathcal{L}_{\T}(X)\coloneq\sum_{e\in\T}\delta_e(X)
\]
for the sum of the truncated lengths of the arcs of $\T$ on $X$. Our main result is an analogue of Wolpert's bound in which the lengths of the curves of a pants decomposition are replaced by the truncated lengths of the arcs of an ideal triangulation.

\begin{theorem*}[\ref{mainthm}]
Let $g,n$ be integers with $n\geq1$ and $2g-2+n>0$. There is a universal constant $C<6$ such that any two hyperbolic surfaces $Y,Z\in\teich_{g,n}$, equipped with the same ideal triangulation $\T$, satisfy
\[
    \dwp(Y,Z)
    \le C\sqrt{n+\mathcal{L}_{\T}(Y)+\mathcal{L}_{\T}(Z)}.
\]
\end{theorem*}

More precisely, $C=\sqrt{8+2\pi B}<5.82$, where $B=2\sqrt2+\log(1+\sqrt2)/\log2$. The constant does not depend on the topology of the surface nor on the geometry of $Y$ and $Z$. A different approach to Theorem~\ref{mainthm} can be based on the formula of Wolf for the Weil--Petersson Hessian of the truncated lengths of geodesic arcs joining cusps~\cite[Section~6.1 and Theorem~6.4]{Wolf2012}.

Shears and truncated lengths are related through Penner's $h$-lengths~\cite{Penner2012}: the shear points of an arc lie on its truncated part, so that $|s_e(X)|\leq\delta_e(X)-2\log4$, and conversely the truncated lengths are controlled by the sums of shears along consecutive arcs around each cusp (Lemma~\ref{lem:shear-trunc}). Let $\mathcal{O}_{\T}(X)$ be the maximum absolute value of a sum of shearing coordinates along consecutive arcs ending in the same cusp (with arc endpoints counted with multiplicity). This is the oscillation norm of Kahn--Markovi\'c~\cite{KahnMarkovic2008}, see Section~\ref{subsec:shear}. As a consequence of the main theorem we obtain the following upper bound 

\begin{cor*}[\ref{mainthm2}]
Let $g,n$ be integers with $n\geq1$ and $2g-2+n>0$. Any two hyperbolic surfaces $Y,Z\in\teich_{g,n}$, equipped with the same ideal triangulation $\T$, satisfy
\[
    \dwp(Y,Z)\leq 21\sqrt{(2g-2+n)\Bigl(\log\bigl(24(2g-2+n)\bigr)+\mathcal{O}_{\T}(Y)+\mathcal{O}_{\T}(Z)\Bigr)}.
\]
\end{cor*}
Kahn--Markovi\'c~\cite[Theorem~2.1]{KahnMarkovic2008} obtained WP distance estimates in terms of oscillation norm and the proportion of edges whose shear exceeds a prescribed threshold, with a constant depending on an upper bound for the individual shears.

Since the sum of all shearing coordinates along arcs ending at the same cusp is zero, we have $\mathcal{O}_{\T}(X)\leq\sum_{e\in\T}|s_e(X)|$. Corollary~\ref{mainthm2} therefore yields

\begin{cor*}[\ref{cor:shear-l1}]
Let $g,n$ be integers with $n\geq1$ and $2g-2+n>0$.
Any two surfaces $Y,Z\in\teich_{g,n}$, equipped with the
same ideal triangulation $\T$, satisfy
\[
    \dwp(Y,Z)
    \leq 21\sqrt{(2g-2+n)\left(\log(24(2g-2+n))
             +\sum_{e\in\T}\bigl(|s_e(Y)|+|s_e(Z)|\bigr)\right)}.
\]
\end{cor*}

\addtocontents{toc}{\SkipTocEntry}
\subsection*{Outline of the proof}
Let $Y$ and $Z$ carry the same marked ideal triangulation $\T$. Truncate the cusps of $Y$ and $Z$ along the horocycles of length $1/4$ and glue the truncated surface to its mirror image, inserting at each cusp a Euclidean cylinder of height $T$ between the two corresponding horocycles. This gives two closed Riemann surfaces $W_T(Y)$ and $W_T(Z)$ of genus $2g+n-1$. Following an arc $e$ of $\T$ on the truncated surface, then its mirror image, and closing up with two vertical segments in the cylinders gives a closed curve $d(e)$. The curves $d(e)$, $e\in\T$, give a pants decomposition $\mathcal{P}_{\T}$ of $W_T(Y)$ and $W_T(Z)$. 

As $T\to\infty$ the core curves of the inserted cylinders are pinched. By using Wolpert's plumbing coordinates~\cite{Wolpert1990}, we show that $W_T(Y)$, resp. $W_T(Z)$, converges in the augmented Teichmüller space of the double to the noded surface $D(Y)$, resp. $D(Z)$, each obtained from the corresponding surface and its mirror image by pairing corresponding cusps into nodes. Denote by $d_{\WP}^{\mathrm{aux}}$ the WP distance on the augmented Teichmüller space of the closed double, of genus $2g+n-1$. The noded surfaces $D(Y)$ and $D(Z)$ lie in the stratum where exactly the core curves of the cylinders are pinched. This stratum is isometric to the product of the Teichmüller spaces of the two halves~\cite{Masur1976,Wolpert2003} so that
\[
    \sqrt2\,\dwp(Y,Z)=d_{\WP}^{\mathrm{aux}}\bigl(D(Y),D(Z)\bigr)\leq d_{\WP}^{\mathrm{aux}}\bigl(D(Y),P_{\T}\bigr)+d_{\WP}^{\mathrm{aux}}\bigl(P_{\T},D(Z)\bigr).
\]
We bound each term by passing through the surface $W_{T_0}(V)$ with $T_0=2\pi$,
\begin{align*}
    \sqrt2\,\dwp(Y,Z)\leq &\, d_{\WP}^{\mathrm{aux}}(D(Y),W_{T_0}(Y)) +  d_{\WP}^{\mathrm{aux}}(W_{T_0}(Y),P_{\T}) \\
    &+ d_{\WP}^{\mathrm{aux}}(P_{\T},W_{T_0}(Z)) + d_{\WP}^{\mathrm{aux}}(W_{T_0}(Z),D(Z)).
\end{align*}
The first term is the cost of opening the cylinders along a plumbing deformation (Lemma~\ref{lem:opening}) while the second term is bounded by Wolpert's estimate applied to the pants decomposition $\mathcal{P}_{\T}$ of $W_{T_0}(Y)$, respectively $W_{T_0}(Z)$ namely $\sqrt{2\pi\sum_{e\in\T}\ell_{W_{T_0}(V)}(d(e))}$. The length of each curve $d(e)$ is controlled by the truncated length of $e$ (Lemma~\ref{lem:arcs}): $\ell_{W_{T_0}(V)}(d(e))\leq2\sqrt2\,\delta_e(V)+4\log(1+\sqrt2).$
\smallbreak
The proof combines classical tools, for instance Wolpert used doubles and opened nodes to prove the strict WP convexity of the distance between horocycles of fixed lengths~\cite[Theorem~6.1]{Wolpert2009}. He also studied shear deformations as limits of opposing twists on opened doubles~\cite[Sections~4--5]{Wolpert2015}.
Mondello gave a formula for the WP Poisson structure in arc-length coordinates on the Teichm\"uller space of compact hyperbolic surfaces with geodesic boundary by doubling the arcs of a triangulation to obtain a pants decomposition~\cite[Theorem~4.1]{Mondello2009}.

\addtocontents{toc}{\SkipTocEntry}
\subsection*{Acknowledgments}
I am grateful to my supervisor Hugo Parlier for helpful discussions and for his encouragement.

\addtocontents{toc}{\SkipTocEntry}
\subsection*{Use of AI tools}
The idea behind this work is to improve the bound on the diameter of the hexagon graph obtained in the first version of~\cite{abadie2026}. My role was that of an archaeologist who knows where to dig: from~\cite{abadie2026}, I identified the precise step needed for this improvement and the part of the literature where it should be found, and I used AI tools as assistants to unearth the argument and to clean and polish it. More precisely, I asked GPT-6-Astra and GPT-6.1-Sol for possible approaches to this step. They suggested a nice strategy, relying on results of Wolpert~\cite{Wolpert1990,Wolpert2003}. I then polished this strategy into the proofs presented here. During this process, I used Claude Opus 5.5 as an assistant to improve the writing, and to help check the clarity and consistency of notation, statements and proofs.

\section{Set up}\label{sec:background}
This section fixes notation and recalls the background material we use. Throughout, $\Sigma$ is an oriented surface of genus $g$ with $n$ punctures and Euler characteristic $\chi(\Sigma)=2-2g-n<0$.

\addtocontents{toc}{\SkipTocEntry}
\subsection{Weil--Petersson distance}\label{subsec:WP} The \emph{Teichmüller space} $\teich_{g,n}=\teich(\Sigma)$ is the space of isotopy classes of complete finite-area hyperbolic metrics on $\Sigma$. Equivalently, it is the space of marked complete finite-area hyperbolic surfaces $(X,f)$, where $f\colon\Sigma\to X$ is a homeomorphism, up to the equivalence $(X,f)\sim(Y,g)$ if and only if $g\circ f^{-1}\colon X\to Y$ is isotopic to an isometry.

Let $X\in\teich_{g,n}$ equipped with the hyperbolic metric $\rho^2|dz|^2$ with area form $dA_X=\rho^2\,dx\,dy$. The cotangent space of $\teich_{g,n}$ at $X$ is the space $Q(X)$ of holomorphic quadratic differentials on $X$.
Let $B(X)$ be the space of Beltrami differentials on $X$ and let $\mu \in B(X)$. Let $\phi \in Q(X)$, there is a bilinear pairing between $B(X)$ and $Q(X)$ given by
\[
    \langle \mu, \phi \rangle = \int_X \mu \phi.
\]
Denote by $N(X) \subset B(X)$ the subset 
\[
    N(X) = \{ \mu \in B(X) \mid \langle \mu, \phi \rangle = 0 \text{ for all } \phi \in Q(X) \}.
\] 
The tangent space of $\teich_{g,n}$ at $X$ is identified with the quotient $B(X)/N(X)$. Every class in $B(X)/N(X)$ contains a unique harmonic Beltrami differential which is of the form $\mu = \rho^{-2} \overline{\psi}$ for some $\psi \in Q(X)$. 
The Weil--Petersson (WP) distance is induced by the norm 
\[
    \left\|[\mu]\right\|_{\WP}=\sup_{\varphi\in Q(X)\smallsetminus\{0\}}\frac{|\langle \mu, \varphi \rangle|}{\|\varphi\|_2},
    \qquad
    \|\varphi\|_2^2=\int_X\frac{|\varphi|^2}{\rho^2}.
\]
Writing $\mu\varphi=(\mu\rho)(\varphi/\rho)$, the Cauchy--Schwarz inequality gives, for every $\mu\in B(X)$,
\begin{equation}\label{eq:beltrami-bound}
    \left\|[\mu]\right\|_{\WP}^2\leq\int_X|\mu|^2\,dA_X,
\end{equation}
with equality when $\mu$ is harmonic, see for instance~\cite[Section 2]{BBB2023}.

\addtocontents{toc}{\SkipTocEntry}
\subsection{The augmented Teichmüller space}\label{subsec:augmentedTeich}The Weil--Petersson metric is not complete~\cite{Wolpert1975,Chu1976}, its completion, the \emph{augmented Teichmüller space} $\overline{\teich}_{g,n}$, is obtained by adding the noded surfaces which arise as limits of hyperbolic surfaces along which the lengths of the curves of some multicurve $\Gamma$ tend to zero~\cite{Bers1974,Masur1976,Abikoff1977,Wolpert2003}:
\[
    \overline{\teich}_{g,n}=\teich_{g,n}\cup\bigcup_{\Gamma\neq\emptyset}S(\Gamma).
\]
The \emph{stratum} $S(\Gamma)$ consists of the noded surfaces at which exactly the curves of $\Gamma$ are pinched. It is the Teichmüller space of the surface $\Sigma\smallsetminus\Gamma$ (possibly disconnected), that is the product $\prod_i\teich(\Sigma_i)$ over the components $\Sigma_i$ of $\Sigma\smallsetminus\Gamma$. When $\Gamma=\Ps$ is a pants decomposition, $S(\Ps)$ is a single \emph{maximally noded surface}. The WP metric extends to $\overline{\teich}_{g,n}$~\cite{Masur1976}, and its restriction to a stratum $S(\Gamma)$ is the product of the Weil--Petersson metrics of the Teichmüller spaces of the components of $\Sigma\smallsetminus\Gamma$~\cite[Theorem~2(1), p.~17]{Wolpert2003}. The WP metric near a stratum $S(\Gamma)$ is close to the product of the Weil--Petersson metrics on its components~\cite{Masur1976, Wolpert2003, DaskalopoulosWentworth2003, Yamada2004}.

A point $p\in\overline{\teich}(\Sigma)$ is represented by a marked noded surface $(X,f)$. If the marking $f$ collapses the curves $\alpha_1,\ldots,\alpha_k$ to the nodes of $X$, define
\[
    \Gamma(p)=\{[\alpha_1],\ldots,[\alpha_k]\},
\]
where $[\alpha_i]$ denotes the free homotopy class of $\alpha_i$ on $\Sigma$. Thus $\Gamma(p)$ is the multicurve on the reference surface specifying the nodes of $p$. By~\cite[Theorem~13, p.~31]{Wolpert2003} any two points $p,q\in\overline{\teich}_{g,n}$ are joined by a unique length-minimizing WP geodesic, and its interior is contained in the stratum $S(\Gamma(p)\cap\Gamma(q))$. Consequently, the ambient WP distance between two points of $S(\Gamma)$ equals their intrinsic distance for the product WP metric. 

The WP distance to a stratum is controlled by the lengths of the pinched curves~\cite{Wolpert2008}, see also~\cite[Corollary 21]{Wolpert2003}. Hence, for every $X\in\teich_{g,n}$ and every multicurve $\Gamma$,
\begin{equation}\label{eq:pinch}
    \dwp(X,S(\Gamma))\leq\sqrt{2\pi\sum_{\gamma\in\Gamma}\ell_X(\gamma)}.
\end{equation}
When $\Gamma=\Ps$ is a pants decomposition, the left-hand side is the distance from $X$ to the maximally noded surface $S(\Ps)$.

The mapping class group $\Mod(\Sigma)$ acts on $\teich_{g,n}$ by changing the marking, and the quotient is the \emph{moduli space} $\Ms_{g,n}=\teich_{g,n}/\Mod(\Sigma)$. The WP metric is invariant under the action of $\Mod(\Sigma)$ and thus descends to $\Ms_{g,n}$. The quotient $\overline{\teich}_{g,n}/\Mod(\Sigma)$ is compact. 

\addtocontents{toc}{\SkipTocEntry}
\subsection{Pants decomposition and Ideal triangulation}\label{subsec:pants-triang}
A \emph{pants decomposition} $\Ps$ of $\Sigma$ is a maximal collection of disjoint, pairwise non-homotopic, essential simple closed curves. It contains $3g-3+n$ curves and $\Sigma\smallsetminus\Ps$ consists of $2g-2+n$ pairs of pants. On a hyperbolic surface we always use geodesic representatives.

Assume that $\Sigma$ has $n\geq 1$ puncture(s). An \emph{ideal triangulation} $\T$ of $\Sigma$ is a maximal collection of disjoint, pairwise non-homotopic, essential simple arcs with both endpoints at punctures. It contains
\[
    E\coloneq6g-6+3n
\]
arcs and $\Sigma\smallsetminus\T$ consists of $4g-4+2n$ triangles. An arc may have both endpoints at the same puncture, and so a triangle may have two sides on the same arc or several vertices at the same puncture. For a puncture $c$, we write $v(c)$ for the number of endpoints of arcs at $c$, an arc with both endpoints at $c$ being counted twice, so that $\sum_c v(c)=2E$.

Truncating the cusps along disjoint horocycles turns each ideal triangle into a \emph{truncated hexagon}: a right-angled hexagon whose sides alternate between three geodesic segments, contained in the arcs, and three horocyclic segments, contained in the corners.

\addtocontents{toc}{\SkipTocEntry}
\subsection{Fenchel--Nielsen coordinates, Shearing coordinates}\label{subsec:shear}
Given a pants decomposition $\Ps$, the \emph{Fenchel--Nielsen coordinates} associate to $X\in\teich_{g,n}$ the lengths and the twist parameters along the curves of $\Ps$,
\[
    \FN_{\Ps}\colon\teich_{g,n}\longrightarrow(\R_{>0}\times\R)^{|\Ps|},\qquad X\longmapsto(\ell_X(\gamma),\tau_X(\gamma))_{\gamma\in\Ps},
\]
and $\FN_{\Ps}$ is a real-analytic homeomorphism. These coordinates extend to the strata $S(\Gamma)$ with $\Gamma\subseteq\Ps$ by allowing curves in $\Gamma$ to have length zero~\cite{Abikoff1980}. When a curve has length zero, there is no twist around it, hence the \emph{extended} length-twist parameters take value in the space
\begin{equation*}
    \frac{\R_{\geq 0}\times\R}{\{0\}\times\R},
\end{equation*}
where the topology around $\Zero\coloneq[\{0\}\times\R]$ is defined as follows: $(\ell_n,\tau_n)\to\Zero$ if and only if $\ell_n\to 0$\footnote{Note that this is \emph{distinct} from the quotient topology on that space. In the quotient topology, we would have $(\ell_n,\tau_n)\to\Zero$ if and only if $\ell_n\to 0$ and $(\tau_n)_{n\in\N}$ is bounded. In the augmented Teichmüller space, the twists can take arbitrary values when $\ell_n\to 0$.}. For a point of $S(\Gamma)$, the \emph{extended Fenchel--Nielsen coordinates} along the curves of $\Ps\smallsetminus\Gamma$ are the Fenchel--Nielsen coordinates of its components with respect to the pants decomposition induced by $\Ps$ on $\Sigma\smallsetminus\Gamma$.

Similarly, given an ideal triangulation we can define coordinates on the Teichmüller space. 
All ideal triangles are isometric to the ideal triangle $(0,1,\infty)$ in $\Hh^2$, and each has an inscribed circle tangent to its three sides. These three tangency points are called the \emph{shear points}. Let $\Delta,\Delta'$ be two ideal triangles sharing a common side $e$. The \emph{shear} parameter $s_e(X)$ is the signed distance along $e$ between the shear points of $\Delta$ and of $\Delta'$ on $e$. 

Let $\T$ be an ideal triangulation of a hyperbolic surface of type $(g,n)$ with $E=6g-6+3n$ ideal arcs. The shear parameters along the arcs of $\T$ are not independent: at each cusp $c$, the sum of the shears along the arcs ending at $c$ is zero, an arc with both endpoints at $c$ being counted twice. The \emph{shearing coordinates} are the shear parameters along all the arcs of $\T$. They give a real-analytic homeomorphism from the Teichmüller space onto the subspace of $\R^{E}$ defined by these $n$ relations, which has dimension $6g-6+2n$~\cite{Bonahon1996,fock1998,Thurston1998,Penner2012},
\[
    \teich_{g,n}\longrightarrow\R^{E},\qquad X\longmapsto(s_e(X))_{e\in\T}.
\]
Following Kahn and Markovi\'c~\cite{KahnMarkovic2008}, for every cusp $c$ we list the endpoints of arcs at $c$ in cyclic order, $e_{c,1},\dots,e_{c,v(c)}$, and we write $s_{e_{c,r}}(X)$ for the shear parameter along the arc containing $e_{c,r}$. The \emph{oscillation norm} of $X$ with respect to $\T$ is
\[
    \mathcal{O}_{\T}(X)\coloneq\max_{c}\;\max_{1\leq j\leq k\leq v(c)}\;\Bigl|\sum_{r=j}^{k}s_{e_{c,r}}(X)\Bigr|.
\]
In particular, the sums of the shears over cyclically consecutive endpoints at $c$ are bounded by $\mathcal{O}_{\T}(X)$ in absolute value.

\addtocontents{toc}{\SkipTocEntry}
\subsection{Cusp truncations}\label{subsec:truncation} 
Each cusp $c$ of $X\in\teich_{g,n}$ has an open neighborhood isometric to
\[
    N_2(c)=\{z=x+iy\in\Hh \mid y>1/2\}/(z\mapsto z+1),
    \qquad ds=\frac{|dz|}{y},
\]
bounded by a horocycle of length $2$ and called \emph{standard cusp neighborhood}. For $y_0> 1/2$, the horocycle at height $y_0$
has length $1/y_0$. Throughout we fix
\[
    h=\frac{1}{4},\qquad a=\frac{1}{h}=4,
\]
and we truncate each cusp along the horocycle $\{y=a\}$, of length $h$. The open cusp neighborhoods $N_h(c)=\{y> a\}$ are contained in the standard ones, and they are embedded and pairwise disjoint~\cite{Keen1974},~\cite[Theorem~4.4.6]{Buser}. We set
\[
    X^h=X\smallsetminus\bigcup_c N_h(c).
\]
Thus $X^h$ is a compact surface with boundary, and
$\partial X^h$ consists of the $n$ horocycles $\{y=a\}$,
each of length $h$.

Let $e$ be an ideal arc of $\T$, realized as a geodesic on $X$. Since $e$ is simple, it meets each cusp neighborhood $N_h(c)$ only along terminal rays, see~\cite[proof of Lemma 8.8]{abadie2026}, so $e\cap X^{h}$ is a single geodesic segment, orthogonal to $\partial X^{h}$. Its length is the \emph{truncated length} $\delta_e(X)$ of $e$, that is
\begin{equation*}
    \delta_e(X) \coloneq\ell_X(e\cap X^{h}).
\end{equation*}
This is the truncated length of~\cite[Definition 2.4]{abadie2026} for the empty multicurve and horocycles of length $h$. In particular $\delta_e$ is a continuous function on $\teich_{g,n}$~\cite[Lemma 3.1]{abadie2026}. We set
\[
    \mathcal{L}_{\T}(X)\coloneq\sum_{e\in\T}\delta_e(X).
\]
Truncated lengths can be expressed in terms of $\lambda$-lengths and $h$-lengths~\cite[Chapter 1, Section 4]{Penner2012}. Let $\h_1,\h_2$ be two horocycles centered at distinct ideal points of $\Hh$, let $e$ be the geodesic joining their centers, and let $\delta$ be the signed distance along $e$ between $\h_1\cap e$ and $\h_2\cap e$, positive if and only if $\h_1$ and $\h_2$ are disjoint. The \emph{$\lambda$-length} of $\h_1,\h_2$ is
\[
    \lambda(\h_1,\h_2)=e^{\delta/2}.
\]
A \emph{decoration} of an ideal triangle $\Delta\subset\Hh$ is a choice of a horocycle centered at each of its three vertices. For a vertex $v$ of $\Delta$, the segment of the horocycle at $v$ contained in $\Delta$ is the \emph{corner} of $\Delta$ at $v$, and its hyperbolic length $h_v(\Delta)$ is the \emph{$h$-length} of this corner. The $h$-lengths of a decorated ideal triangle are determined by the $\lambda$-lengths of its sides, and conversely.

\begin{lemma}[{\cite[Chapter 1, Lemma 4.9]{Penner2012}}]\label{lem:corners}
    Let $\Delta$ be a decorated ideal triangle and let $e$ be a side of $\Delta$ with endpoints $u$ and $w$. If $\lambda_e$ denotes the $\lambda$-length of the horocycles at $u$ and $w$, then $h_u(\Delta)\,h_w(\Delta)=\lambda_e^{-2}$. In particular, the signed length $\delta_e$ of the segment of $e$ between the horocycles at $u$ and $w$ satisfies
    \begin{equation}\label{eq:corner-identity}
        \delta_e=\log\frac{1}{h_u(\Delta)}+\log\frac{1}{h_w(\Delta)}.
    \end{equation}
\end{lemma}

When $\Delta$ is a lift of a triangle of $\T$ on $X$, decorated by the lifts of the horocycles of length $h$, its corners project to corners of the triangles of $\T$. In particular every corner has $h$-length at most $h$, and Lemma~\ref{lem:corners} gives, for every $X\in\teich_{g,n}$ and every arc $e\in\T$,
\begin{equation}\label{eq:delta-lower}
    \delta_e(X)\geq2\log\frac1h=2\log4.
\end{equation}

The following lemma relates shear parameters and truncated lengths.

\begin{lemma}\label{lem:shear-trunc}
    Let $X\in\teich_{g,n}$ and let $c$ be a puncture with $v=v(c)$. For $k=1,\dots,v$, let $\Delta_k$ be the triangle whose corner at $c$ lies between the endpoints $e_{c,k}$ and $e_{c,k+1}$, indices modulo $v$, and set $\sigma_k=s_{e_{c,1}}(X)+\dots+s_{e_{c,k}}(X)$. Then, for one of the two orientations of the cyclic order at $c$,
    \begin{equation}\label{eq:corner-shear}
        h_c(\Delta_k)=h\,\frac{e^{\sigma_k}}{\sum_{j=1}^{v}e^{\sigma_j}}.
    \end{equation}
    Consequently, every arc $e\in\T$, with endpoints at the punctures $c$ and $c'$, satisfies
    \[
        |s_e(X)|\leq\delta_e(X)-2\log4
        \qquad\text{and}\qquad
        \delta_e(X)\leq2\log4+\log\bigl(v(c)v(c')\bigr)+2\,\mathcal{O}_{\T}(X),
    \]
    and in particular $\mathcal{L}_{\T}(X)\leq2E\bigl(\log(8E)+\mathcal{O}_{\T}(X)\bigr)$.
\end{lemma}

\begin{proof}
    Lift $c$ to $\infty$ in the upper half-plane $\Hh^2$ such that the stabilizer of $\infty$ is generated by $z\mapsto z+1$. The lifts of the arcs ending at $c$ are vertical lines $\{x=x_k\}$ with $x_1< \dots <x_v<x_1+1$ and such that the intersection point $e_{c,k}$ lifts to the intersection of the vertical line $\{x=x_k\}$ with the horocyle $\{y=a\}$. Then for $k=1,\dots, v$ the triangle $\Delta_k$ lifts to the triangle $(x_k,x_{k+1},\infty)$ where $x_{v+1}=x_1+1$. Set $\lambda_k=x_{k+1}-x_k$, so that $\lambda_1+\dots+\lambda_v=1$, and $\lambda_0=\lambda_v$. 
    \begin{figure}[h]
    \centering
    \begin{overpic}[width=1.0\linewidth,keepaspectratio]{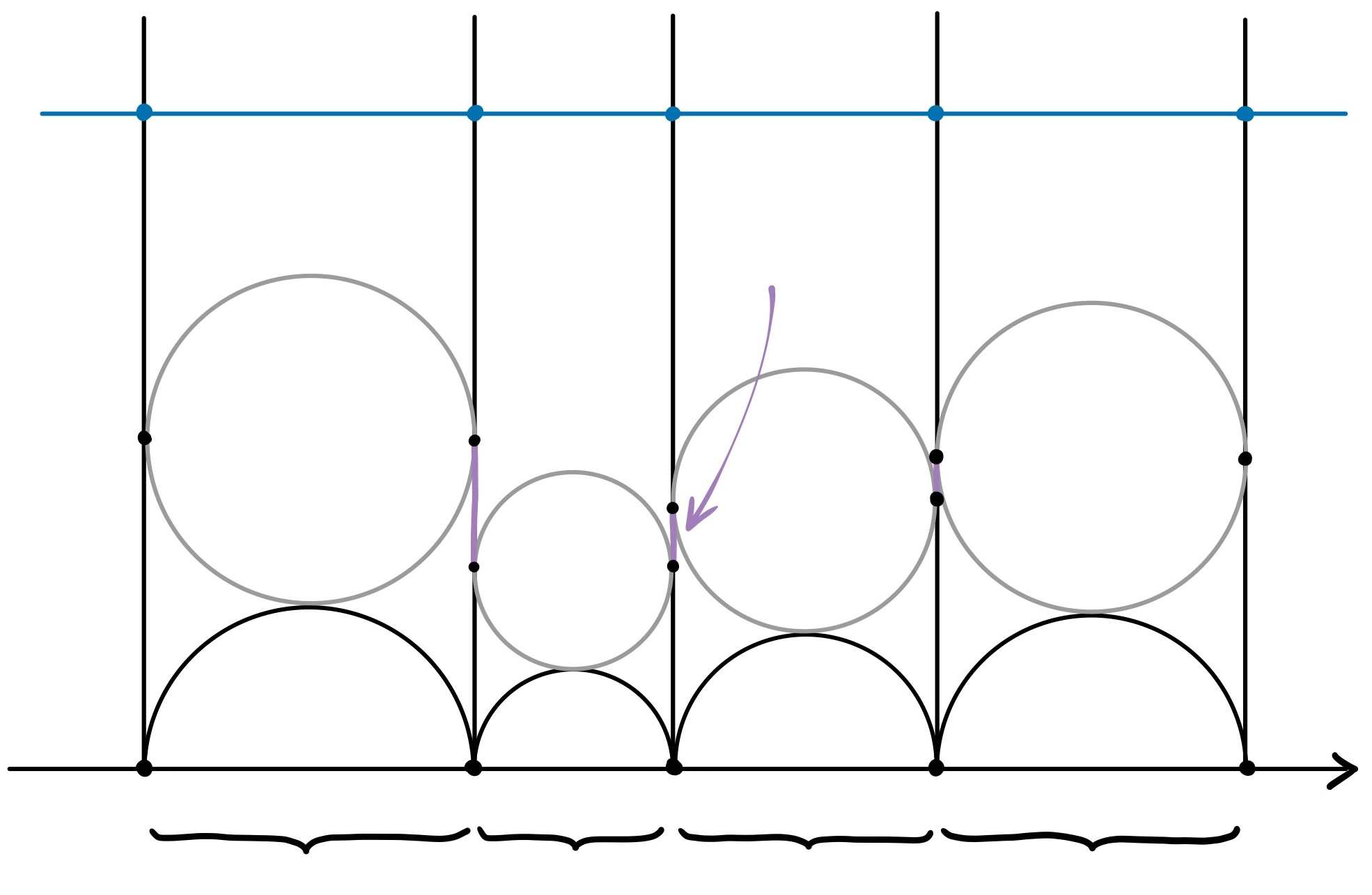}
        \put(21,33){$\Delta_1$}
        \put(41,23){$\Delta_2$}
        \put(58,28){$\Delta_3$}
        \put(78,31){$\Delta_4$}
        \put(11,55){\color{MidnightBlue}$e_{c,1}$}
        \put(35,55){\color{MidnightBlue}$e_{c,2}$}
        \put(49.5,55){\color{MidnightBlue}$e_{c,3}$}
        \put(69,55){\color{MidnightBlue}$e_{c,4}$}
        \put(95,58){\color{MidnightBlue}$y=a$}
        \put(20,58){\color{MidnightBlue}$h\lambda_1$}
        \put(40,58){\color{MidnightBlue}$h\lambda_2$}
        \put(57,58){\color{MidnightBlue}$h\lambda_3$}
        \put(77,58){\color{MidnightBlue}$h_c(\Delta_4)=h\lambda_4$}
        \put(10,7){$x_1$}
        \put(33,7){$x_2$}
        \put(48,7){$x_3$}
        \put(67,7){$x_4$}
        \put(87,7){$x_1+1$}
        \put(21,1){$\lambda_1$}
        \put(41,1){$\lambda_2$}
        \put(58,1){$\lambda_3$}
        \put(78,1){$\lambda_4$}
        \put(49.25,45){\color{Orchid}$s_{e_{c,3}}=\log\left(\frac{\lambda_3}{\lambda_2}\right)$}
        \end{overpic}\vspace*{-0.25cm}
        \caption{A fan of edges ending in the cusp $c$ lifted to $\infty$ in the proof of Lemma~\ref{lem:shear-trunc} for $v(c)=4$.}\label{Fig-fan-edges}
    \end{figure}
    The map $z\mapsto\lambda_k z+x_k$ sends the triangle $(0,1,\infty)$ to the lift of $\Delta_k$. See Figure~\ref{Fig-fan-edges}. Since the shear points of the former on its vertical sides are $i$ and $1+i$, the shear points of $\Delta_k$ on its vertical sides are $x_k+i\lambda_k$ and $x_{k+1}+i\lambda_k$. Therefore the two shear points on $\{x=x_k\}$ are at heights $\lambda_{k-1}$ and $\lambda_k$, and for one of the two orientations of the cyclic order 
    \[
        s_{e_{c,k}}(X)=\log(\lambda_k/\lambda_{k-1}).
    \] 
    Hence $\sigma_k=s_{e_{c,1}}(X)+\dots+s_{e_{c,k}}(X)=\log(\lambda_1/\lambda_{0})+\dots \log(\lambda_k/\lambda_{k-1})=\log(\lambda_k/\lambda_{v})$ and so 
    \begin{equation}\label{eq:lambda-SumShear}
        \lambda_k=\lambda_v e^{\sigma_k}.
    \end{equation}
    Consequently $1=\sum_{k=1}^v \lambda_k=\sum_{k=1}^v \lambda_ve^{\sigma_k}$ which yields 
    \begin{equation}\label{eq:lambdav-sum-sigma}
        \lambda_v=\frac{1}{\sum_j e^{\sigma_j}}.
    \end{equation}
    The corner of $\Delta_k$ at $c$ is the horizontal segment $\{x+ia : x_k\leq x\leq x_{k+1}\}$, of Euclidean length $\lambda_k$ at height $a$ so we obtain
    \begin{equation}\label{eq:h-lambda}
        h_c(\Delta_k)=\int_{x_k}^{x_{k+1}}\frac{dx}{a}=\frac{\lambda_k}{a}=h\lambda_k,
    \end{equation}
    and therefore
    \[
        h_c(\Delta_k)\underset{\eqref{eq:h-lambda}}{=}h\lambda_k\underset{\eqref{eq:lambda-SumShear}}{=}h\cdot\lambda_v e^{\sigma_k}\underset{\eqref{eq:lambdav-sum-sigma}}{=}h\,\frac{e^{\sigma_k}}{\sum_{j=1}^{v}e^{\sigma_j}}.
    \]
    Let $e\in\T$ and choose an endpoint of $e$ at $c$, say $e_{c,k}$, so that $\{x=x_k\}$ is a lift of $e$. On $\{x=x_k\}$, the horocycle of length $h$ at $c$ lies at height $a$, while the shear points of $\Delta_{k-1}$ and $\Delta_k$ lie at heights $\lambda_{k-1},\lambda_k$. Since $\lambda_{k-1},\lambda_k\leq\lambda_1+\dots+\lambda_v=1<a$, both shear points lie below this horocycle, that is outside $N_h(c)$, and their distances to it are
    \[
        \log\frac{a}{\lambda_{k-1}}\geq\log a=\log4
        \qquad\text{and}\qquad
        \log\frac{a}{\lambda_k}\geq\log4.
    \]
    Applying the same argument at the other endpoint of $e$, at a puncture $c'$ (possibly $c'=c$), shows that they also lie outside $N_h(c')$, at distance at least $\log4$ from its boundary. Hence both shear points lie on the truncated segment $e\cap X^h$, of length $\delta_e(X)$, at distance at least $\log4$ from each of its endpoints. Since $|s_e(X)|$ is the distance between them, $|s_e(X)|\leq\delta_e(X)-2\log4$.

    We now prove the upper bound on $\delta_e(X)$. By~\eqref{eq:corner-shear},
    \[
        \log\frac{1}{h_c(\Delta_k)}=\log\frac{1}{h}+\log\sum_{j=1}^{v}e^{\sigma_j-\sigma_k},
    \]
    the sum contains the term $e^0=1$, and each $\sigma_j-\sigma_k$ is, up to sign, a sum of shears over consecutive arcs ending at $c$, so $|\sigma_j-\sigma_k|\leq\mathcal{O}_{\T}(X)$. Hence
    \[
        \log\frac1h\leq\log\frac{1}{h_c(\Delta_k)} \leq \log\frac1h+\log\left(v e^{\mathcal{O}_{\T}(X)}\right) \leq \log\frac1h+\log v(c)+\mathcal{O}_{\T}(X).
    \]
    Applying this at both endpoints of $e$, in one of its adjacent triangles, together with~\eqref{eq:corner-identity}, yields 
    \[ 
        \delta_e(X)\leq2\log4+\log\bigl(v(c)v(c')\bigr)+2\,\mathcal{O}_{\T}(X).
    \] 
    Finally, since $\sum_{e}\log\bigl(v(c)v(c')\bigr)=\sum_cv(c)\log v(c)\leq2E\log(2E)$ by summing over the arcs of $\T$ we obtain \[\mathcal{L}_{\T}(X)\leq2E\log4+2E\log(2E)+2E\,\mathcal{O}_{\T}(X),\] where $E=6g-6+3n$ is the number of arcs in $\T$.
\end{proof}

\addtocontents{toc}{\SkipTocEntry}
\subsection{Plumbing and degenerating families}\label{subsec:plumbing}
We recall the plumbing construction of degenerating families of Riemann surfaces from~\cite[Section 2]{Wolpert1990}, see also~\cite[Section 3]{Wolpert2003} and~\cite[Section 2]{ObitsuWolpert2008}.
\subsubsection{Surfaces with nodes and cusp coordinates}
A \emph{Riemann surface with nodes} $R$ is a connected one-dimensional complex space in which every point has a neighborhood isomorphic either to the unit disc of $\C$, or to $\{(z,w)\in\C^2\mid zw=0,\; |z|,|w|<1\}$, by an isomorphism sending the point to the origin. The points of the second type are the \emph{nodes} of $R$. The components of $R\smallsetminus\{\text{nodes}\}$ are the \emph{parts} of $R$, and removing a node creates two punctures, which are said to be \emph{paired}. We always assume that the parts have negative Euler characteristic, so that each of them carries a unique complete hyperbolic metric of finite area compatible with its complex structure, in which its punctures are cusps.

Let $p$ be a puncture of a hyperbolic surface $X$. Recall that the standard cusp neighborhood of $p$ and its hyperbolic metric are
\[
    N_2(p)=\{z=x+iy \mid x\in\R/\Z,\ y>1/2\},\qquad \frac{|dz|}{y},
\]
the puncture corresponding to $\infty$, say that $z=x+iy$ is the \emph{normalized coordinate} at $p$. The \emph{cusp coordinate} at $p$ is
\[
    q=e^{2\pi iz}=e^{2\pi ix}\,e^{-2\pi y},
\]
it is a biholomorphism from $N_2(p)$ onto the punctured disc $\{q \in \C \mid 0<|q|<e^{-\pi}\}$ and the puncture corresponds to the center $q=0$ of this disc. By adding the point $p$ to $X$ with $q(p)=0$, which we call \emph{filling the puncture}, $q$ becomes a holomorphic coordinate on the disc $\{q \in \C \mid |q|<e^{-\pi}\}$. In this coordinate, since $|q|=e^{-2\pi y}$, the horocycle $\{y=y_0\}$, of length $1/y_0$, is the circle $\{|q|=e^{-2\pi y_0}\}$. In particular the cusp neighborhood $N_h(p)=\{y> a\}$ removed in the truncation $X^h$ is the punctured disc $\{q \in \C \mid 0<|q|\leq e^{-2\pi a}\}$. Since $dq=(2\pi i)q dz$, we have $|dq|=2\pi|q|\,|dz|$. Moreover $|q|=e^{-2\pi y}$, so $\log|q|=-2\pi y<0$ on $N_2(p)$ and $\bigl|\log|q|\bigr|=2\pi y$. Therefore
    \[
        \frac{|dz|}{y}
        =\frac{|dq|}{2\pi|q|}
        \frac{2\pi}{\bigl|\log|q|\bigr|}
        =\frac{|dq|}{|q|\,\bigl|\log|q|\bigr|}.
    \]
that is, the hyperbolic metric of $N_2(p)$ is the complete hyperbolic metric of the punctured unit disc $\{q \in \C \mid 0<|q|<1\}$ restricted to $\{q \in \C \mid 0<|q|<e^{-\pi}\}$.

The normalized coordinate is unique up to a translation $x\mapsto x+\theta$, $\theta\in\R/\Z$, hence the cusp coordinate is unique up to a rotation $q\mapsto e^{2\pi i\theta}q$. These are the rotationally symmetric coordinates at a puncture of~\cite[Definition 2.1]{Wolpert1990}.

\subsubsection{Discs and annuli around paired punctures} Let $R$ be a surface with nodes $\nu_1,\dots,\nu_k$, and let $R_0=R\smallsetminus\{\nu_1,\dots,\nu_k\}$. Thus $R_0$ is a finite disjoint union of complete finite-area hyperbolic surfaces, and removing the node $\nu_j$ creates a pair of punctures $(p_j,p'_j)$ of $R_0$. Let $U_j\subset R_0$ and $U'_j\subset R_0$ be the standard cusp neighborhoods of $p_j$ and $p'_j$, and let
\[
    q_j\colon U_j\longrightarrow\{w\in\C\mid 0<|w|<e^{-\pi}\},\qquad q'_j\colon U'_j\longrightarrow\{w\in\C\mid 0<|w|<e^{-\pi}\}
\]
be the cusp coordinates, which are biholomorphisms. For $0<r_1<r_2< e^{-\pi}$ we denote by
\[
    D_j(r_2)=\{\xi\in U_j\mid |q_j(\xi)|<r_2\},\qquad \overline D_j(r_1)=\{\xi\in U_j\mid |q_j(\xi)|\leq r_1\},
\]
\[A_j(r_1,r_2)=\{\xi\in U_j\mid r_1<|q_j(\xi)|<r_2\}\]
the corresponding punctured discs and annulus of $R_0$ around $p_j$, and by $D'_j(r_2)$, $\overline D'_j(r_1)$ and $A'_j(r_1,r_2)$ the analogous subsets of $U'_j$ defined with $q'_j$. Fix $0<c< e^{-\pi}$. The punctured discs $D_j(c)$ and $D'_j(c)$, $1\leq j\leq k$, are pairwise disjoint, since they are contained in the standard cusp neighborhoods, see \cite[Theorem 4.4.6]{Buser}. Their complement
\[
    K_c=R_0\smallsetminus\bigcup_{j=1}^k\bigl(D_j(c)\cup D'_j(c)\bigr)
\]
is a surface with boundary consisting of the $2k$ circles
\[
    \{\xi\in U_j\mid |q_j(\xi)|=c\},
    \qquad
    \{\xi'\in U'_j\mid |q'_j(\xi')|=c\},
    \qquad 1\leq j\leq k.
\]
It is compact when all punctures are paired.

\subsubsection{The plumbing construction} Let $\mathbf{t}=(t_1,\dots,t_k)\in\C^k$ with $0<|t_j|<c^2$ for every $j$, so that $|t_j|/c<c$. For each $j$, we write
\[
    A_j=A_j(|t_j|/c,c),\qquad A'_j=A'_j(|t_j|/c,c),\qquad \mathbb A_j=\{w\in\C\mid |t_j|/c<|w|<c\}.
\]
Thus $A_j$ and $A'_j$ are the annuli around the paired punctures $p_j$ and $p'_j$ defined above, and $\mathbb A_j$ is the corresponding annulus in the plane. The cusp coordinates restrict to biholomorphisms $q_j\colon A_j\to\mathbb A_j$ and $q'_j\colon A'_j\to\mathbb A_j$. The map $\iota_j(w)=t_j/w$ is holomorphic on $\C\smallsetminus\{0\}$ and satisfies $|\iota_j(w)|=|t_j|/|w|$, hence maps $\mathbb A_j$ into itself. Moreover $\iota_j\circ\iota_j=\mathrm{id}$, thus $\iota_j$ is a biholomorphism of $\mathbb A_j$, which exchanges its two boundary circles $|w|=c$ and $|w|=|t_j|/c$. Therefore
\[
    \varphi_j=(q'_j)^{-1}\circ\iota_j\circ q_j\colon A_j\longrightarrow A'_j
\]
is a biholomorphism, and for $\xi\in A_j$ and $\xi'\in A'_j$ we have $\xi'=\varphi_j(\xi)$ if and only if $q_j(\xi)\,q'_j(\xi')=t_j$.

Remove from $R_0$ the closed punctured discs $\overline D_j(|t_j|/c)$ and $\overline D'_j(|t_j|/c)$, $1\leq j\leq k$, to obtain the open surface
\[
    R_{\mathbf{t}}^*=R_0\smallsetminus\bigcup_{j=1}^k\Bigl(\overline D_j(|t_j|/c)\cup\overline D'_j(|t_j|/c)\Bigr),
\]
which contains the pairwise disjoint annuli $A_j$ and $A'_j$, $1\leq j\leq k$. For each $j$, we identify every point $\xi\in A_j$ with $\varphi_j(\xi)\in A'_j$, or equivalently, $\xi\sim\xi'$ if and only if $q_j(\xi)\,q'_j(\xi')=t_j$. The \emph{plumbing} of $R_0$ with parameter $\mathbf{t}$ is the quotient $R_{\mathbf{t}}=R_{\mathbf{t}}^*/\!\sim$, endowed with the quotient topology, and we denote by $\pi\colon R_{\mathbf{t}}^*\to R_{\mathbf{t}}$ the quotient map. Since the maps $\varphi_j$ are biholomorphisms, the complex structure of $R_{\mathbf{t}}^*$ descends to $R_{\mathbf{t}}$, which is therefore a Riemann surface, it is closed when every puncture of $R_0$ is paired.

\begin{figure}
    \centering
    \begin{overpic}[width=0.8\linewidth,keepaspectratio]{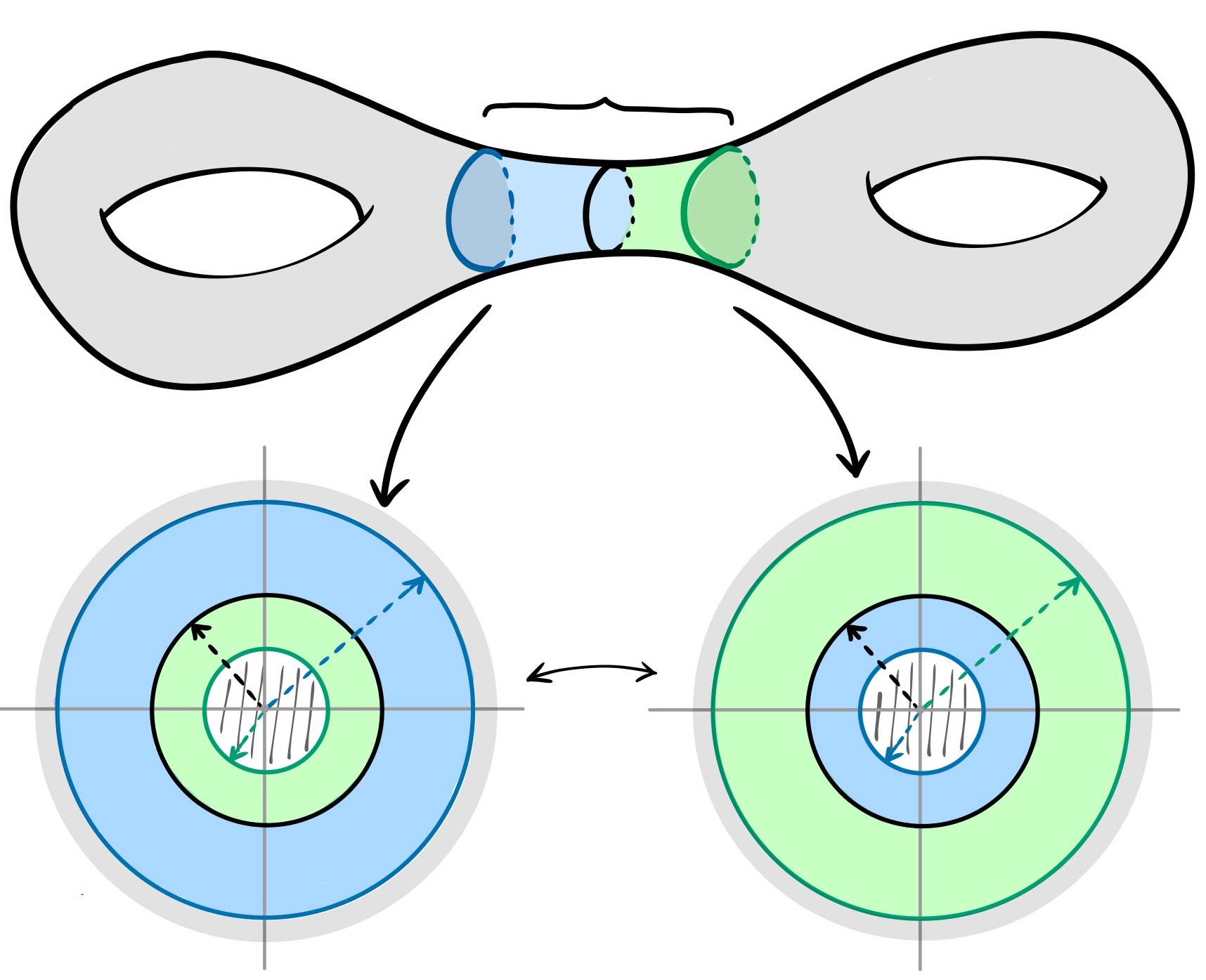}
        \put(3,75){$R_{\mathbf{t}}$}
        \put(49,74){$C_j$}
        \put(11,54){\color{gray}$K_c$}
        \put(81,56){\color{gray}$K_c$}
        \put(33.5,62){\color{MidnightBlue}$\beta_j$}
        \put(57.5,62){\color{ForestGreen}$\beta'_j$}
        \put(43,28){$q_j q_j'=t_j$}
        \put(10,-2){\color{MidnightBlue}$\beta_j=\{|q_j|=c\}$}
        \put(65,-2){\color{ForestGreen}$\beta'_j=\{|q'_j|=c\}$}
        \put(35.5,33.5){\color{MidnightBlue}$c$}
        \put(90,33.5){\color{ForestGreen}$c$}
        \put(15,16){\tiny\color{ForestGreen}$|t_j|/c$}
        \put(70,16){\tiny\color{MidnightBlue}$|t_j|/c$}
        \put(11.5,31){\small$|t_j|^{1/2}$}
        \put(66,31){\small$|t_j|^{1/2}$}
        \end{overpic}
        \caption{The plumbing collar $C_j$, adapted from a figure of~\cite{Matheus2018}.}\label{Fig-Plumbing}
\end{figure}

The $j$-th \emph{plumbing collar} is $C_j=\pi(A_j)=\pi(A'_j)$, see Figure~\ref{Fig-Plumbing}. Every point $x\in C_j$ has exactly two preimages under $\pi$, a point $\xi\in A_j$ and the point $\xi'=\varphi_j(\xi)\in A'_j$. Hence the cusp coordinates induce two holomorphic coordinates on $C_j$, $q_j(x)\coloneq q_j(\xi)$ and $q'_j(x)\coloneq q'_j(\xi')$, with values in $\mathbb A_j$, and by the gluing rule they satisfy $q_j(x)\,q'_j(x)=t_j$. The surface $R_{\mathbf{t}}$ does not depend on the choice of $c$, as long as $|t_j|<c^2$ for every $j$, see~\cite[Remarks in Section~2.4]{Wolpert1990}.

Let $\beta_j=\{\xi\in U_j\mid |q_j(\xi)|=c\}$ and $\beta'_j=\{\xi'\in U'_j\mid |q'_j(\xi')|=c\}$ be the boundary circles of $K_c$. Only points of the open annuli $A_j$ and $A'_j$ are identified, so $K_c$ embeds in $R_{\mathbf{t}}$ and does not meet the collars $C_j$. Note that $|q'_j(\varphi_j(\xi))|=|t_j|/|q_j(\xi)|$ and let $(\xi_n)$ be a sequence in $A_j$. If $|q_j(\xi_n)|\to c$, then, since $q_j$ is a homeomorphism onto its image, $(\xi_n)$ accumulates on $\beta_j$. If $|q_j(\xi_n)|\to|t_j|/c$, then $(\xi_n)$ has no accumulation point in $R_{\mathbf t}^*$ on the side of $A_j$, since the circle $\{|q_j|=|t_j|/c\}$ has been removed. But in $R_{\mathbf t}$ the points $\xi_n$ are identified with $\varphi_j(\xi_n)\in A'_j$, and $|q'_j(\varphi_j(\xi_n))|\to c$, so the sequence accumulates on $\beta'_j$. Hence the closure $\overline C_j$ of $C_j$ in $R_{\mathbf{t}}$ is a closed annulus with boundary $\beta_j\cup\beta'_j$, and
\[
    R_{\mathbf{t}}=K_c\cup\overline C_1\cup\dots\cup\overline C_k,\qquad K_c\cap\overline C_j=\beta_j\cup\beta'_j.
\]
In other words, $R_{\mathbf{t}}$ is obtained from $K_c$ by gluing, for each $j$, a closed annulus along $\beta_j$ and $\beta'_j$.

\subsubsection{The local model} The \emph{plumbing variety} is the complex surface
\[
    V=\bigl\{(u,v,t)\in\C^3\;\big|\;uv=t,\ |u|<1,\ |v|<1,\ |t|<1\bigr\},
\]
with the holomorphic projection $\Pi\colon V\to\{t\in\C\mid|t|<1\}$, $\Pi(u,v,t)=t$. Its \emph{fiber} over $t$ is
\[
    V_t=\bigl\{(u,v)\in\C^2\;\big|\;uv=t,\ |u|<1,\ |v|<1\bigr\}.
\]
For $t\neq0$, the projection $(u,v)\mapsto u$ identifies $V_t$ with the annulus $\{u\in\C\mid |t|<|u|<1\}$. For $t=0$, the fiber
\[
    V_0=\bigl\{(u,0)\in\C^2\;\big|\;|u|<1\bigr\}\cup\bigl\{(0,v)\in\C^2\;\big|\;|v|<1\bigr\}
\]
is the union of two discs meeting at the origin: this is the local model of a node. Thus $V$ is a family of annuli degenerating to a node. It is a local model for $R_{\mathbf{t}}$, near the $j$-th collar: the map
\[
    C_j\longrightarrow V_{t_j},\qquad \xi\longmapsto\bigl(q_j(\xi),q'_j(\xi)\bigr),
\]
is a biholomorphism from the collar $C_j\subset R_{\mathbf{t}}$ onto $\bigl\{(u,v)\in V_{t_j}\;\big|\;|u|<c,\ |v|<c\bigr\}$ which is identified with the annulus $\{u\in\C\mid |t_j|/c<|u|<c\}$ by the projection $(u,v)\mapsto u$.

\subsubsection{The model metric of a plumbing collar} By~\cite[Section 0.1]{Wolpert1990}, for $t\neq0$, the fiber $V_t$, identified with the annulus $\{u\in\C\mid |t|<|u|<1\}$, carries the complete hyperbolic metric
\[
    ds_t=\frac{\pi}{\bigl|\log|t|\bigr|}\csc\left(\frac{\pi\log|u|}{\log|t|}\right)\left|\frac{du}{u}\right|,
\]
and each of the two punctured discs $\{u\in\C\mid 0<|u|<1\}$ and $\{v\in\C\mid 0<|v|<1\}$ of $V_0\smallsetminus\{(0,0)\}$ carries the cusp metric $ds_0=|du|/(|u|\,|\log|u||)$, respectively $|dv|/(|v|\,|\log|v||)$. Writing $u=e^{2\pi i(x+iy)}$, the annulus $\{u\in\C\mid |t|<|u|<1\}$ becomes the cylinder $\{(x,y)\in(\R/\Z)\times\R\mid 0<y<L\}$ with $L=\frac{1}{2\pi}\log\frac{1}{|t|}$, and $ds_t=\lambda_L(y)\,|dz|$, where for $L>0$ we write
\begin{equation}\label{eq:annulus-density}
    \lambda_L(y)=\frac{\pi}{L}\csc\left(\frac{\pi y}{L}\right)=\frac{1}{y}\cdot\frac{\theta}{\sin\theta},\qquad\theta=\frac{\pi y}{L},\quad 0<y<L.
\end{equation}
This is the complete hyperbolic metric of the flat cylinder $\{(x,y)\in(\R/\Z)\times\R\mid 0<y<L\}$.  Hence, for $y$ in compact subsets of $(0,\infty)$, the model metrics $ds_t$ converge to the cusp metric $ds_0=|dz|/y$ as $t\to0$. Under the identification of $C_j$ with part of $V_{t_j}$, the restriction of $ds_{t_j}$ to the collar is the \emph{model metric} of the $j$-th plumbing collar.

Wolpert approximates the hyperbolic metric of $R_{\mathbf{t}}$ by interpolating between the model metrics in the plumbing collars and the hyperbolic metric of $R_0$ away from them~\cite[Sections 3.3--3.4]{Wolpert1990}. The curvature of the resulting \emph{grafted metric} is $-1+O\bigl(\sum_j(\log|t_j|)^{-2}\bigr)$, and a maximum principle compares it with the hyperbolic of $R_{\mathbf{t}}$~\cite[Lemma 3.5]{Wolpert1990}. Expansions of the hyperbolic metric of $R_{\mathbf{t}}$ in the parameters $t_j$ are given in~\cite[Expansion 4.2]{Wolpert1990} and~\cite[Theorem 4]{ObitsuWolpert2008}.

\subsubsection{The plumbing family and plumbing coordinates} The surfaces $R_{\mathbf{t}}$, for $\mathbf{t}$ in the polydisc 
\[
    \Delta_{c^2}^k=\bigl\{\mathbf{t}=(t_1,\dots,t_k)\in\C^k\;\big|\;|t_j|<c^2\ \text{for}\ 1\leq j\leq k\bigr\}
\] 
form a holomorphic family which degenerates to the surface with nodes $R$ at $\mathbf{t}=0$, see~\cite[Sections 0.1 and 2.4]{Wolpert1990}. 
One may also deform the parts: let $\mu_1,\dots,\mu_m$ be Beltrami differentials on $R_0$, supported in $K_c$, whose classes form a basis of the tangent space at $R_0$ of the Teichmüller space of the parts, which has complex dimension $m=3g-3+n-k$, and for $\mathbf{s}\in\C^m$ small let $\mu(\mathbf s)=\sum_{i=1}^ms_i\mu_i$.
Plumbing the deformed parts $R_0^{\mu(\mathbf s)}$ with parameter $\mathbf{t}$ gives surfaces $R_{\mathbf s,\mathbf{t}}$, and $(\mathbf s,\mathbf{t})$ are local coordinates of the compactified moduli space near $R$~\cite[Section 2.4.M]{Wolpert1990},~\cite[Section 3]{Wolpert2003}. 
In this work the parts are fixed, that is $\mathbf s=0$, and we only use the following consequence.

\begin{proposition}[{\cite[Section~3]{Wolpert2003}}]
\label{prop:plumbing}
    Let $R$ be a marked surface with nodes $\nu_1,\dots,\nu_k$, let $R_0=R\smallsetminus\{\nu_1,\dots,\nu_k\}$, and fix $0<c<e^{-\pi}$. For $\mathbf t\in(0,c^2)^k$, let $R_{\mathbf t}$ be the plumbing of $R_0$ with parameter $\mathbf t$, with a marking whose restriction to $K_c$ agrees with the marking of the parts of $R$. Then $R_{\mathbf t}\to R$ in augmented Teichmüller space as $\mathbf t\to0$.
\end{proposition}
\begin{proof}
    By construction, $R$ is the plumbing of the undeformed parts with parameter $\mathbf t=0$, that is, $R=R_{0,0}$ in the plumbing coordinates $(\mathbf{s},\mathbf{t})$.

    By~\cite[Section~3, p.~15]{Wolpert2003}, near $R$ the marked plumbing construction is described by $(\mathbf{s},|\mathbf{t}|,\arg\mathbf{t})$, where each argument $\arg t_j$ is lifted to $\R$. Changing the lift of $\arg t_j$ by $2\pi$ changes the marking by a Dehn twist about the core curve of the collar $C_j$, while $|t_j|\to0$ if and only if the length of this curve tends to zero.

    Recall from Section~\ref{subsec:shear} that, in extended Fenchel--Nielsen coordinates, $(\ell_n,\tau_n)\to\Zero$ if and only if $\ell_n\to0$, with no condition on the twists $\tau_n$. In the same way, $R_{\mathbf{s},\mathbf{t}}\to R$ if and only if $\mathbf{s}\to0$ and $|\mathbf{t}|\to0$, with no condition on the lifted arguments~\cite[Section~3, pp.~11--12]{Wolpert2003}.

    Now let $\mathbf{t}\in(0,c^2)^k$ be sufficiently small. The surface $R_{\mathbf{t}}$ is a plumbing of the undeformed parts, and its marking agrees with that of $R$ on $K_c$. Hence $\mathbf s=0$. As $\mathbf t\to0$, we have $|\mathbf t|\to0$, and any twists in the plumbing collars do not affect convergence to the noded surface. Therefore $R_{\mathbf t}\to R$.
\end{proof}

\section{Doubling and the Weil--Petersson distance}

\subsection{Opening the conformal double}\label{subsec:opening}
Let $Y\in\teich(\Sigma)$ where $\Sigma$ is a surface of type $(g,n)$ with $n\geq 1$. Let $Y^{h}$ be the truncated surface obtained by removing the cusp neighborhoods bounded by the horocycles of length $h$, denote by $z=x+iy$ the normalized coordinate at a cusp $c$. The \emph{mirror} $\overline{Y}$ is the surface $Y$ with the conjugate complex structure. It has the same hyperbolic metric, and at the cusp corresponding to $c$ its normalized coordinate is $z'=-\bar{z}=-x+iy=x'+iy'$.

\subsubsection{Definition of the double}
\begin{figure}[h]
    \centering
    \begin{overpic}[width=0.6\linewidth,keepaspectratio]{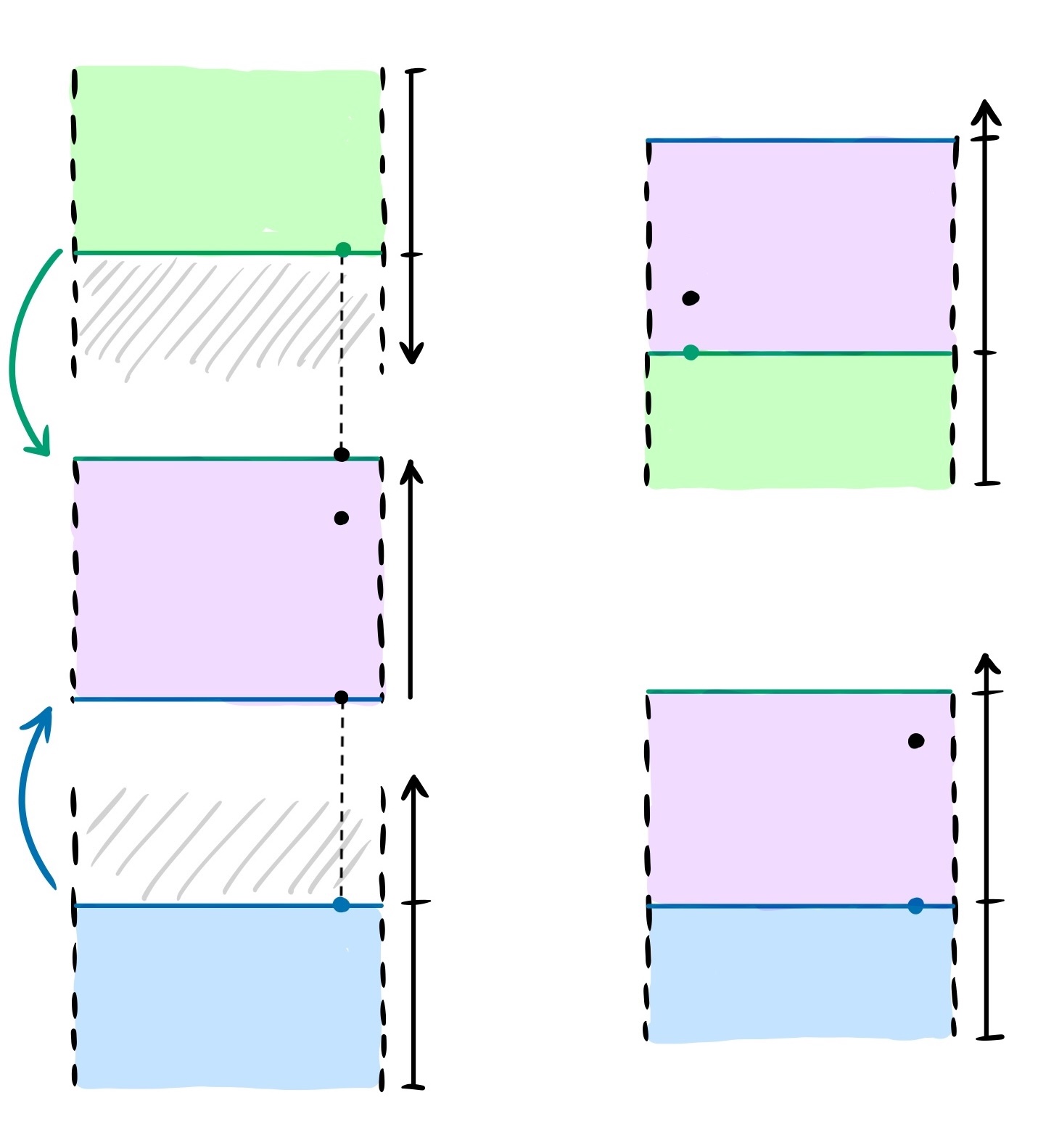}
        \put(-16,73){\tiny\color{ForestGreen}glue $\{s=T\}$}
        \put(-15,70){\tiny\color{ForestGreen}to $\{y'=a\}$}
        \put(-16,33){\tiny\color{MidnightBlue}glue $\{s=0\}$}
        \put(-15,30){\tiny\color{MidnightBlue}to $\{y=a\}$}
        \put(10.8,85.9){\color{ForestGreen}$\overline{Y}^h$}
        \put(25,89){\color{ForestGreen}$\fbox{z'}$}
        \put(10.5,72){\tiny removed$\{y'> a\}$}
        \put(11.5,28){\tiny removed$\{y> a\}$}
        \put(10.8,47.9){\color{Orchid}$C_T$}
        \put(10.8,9.9){\color{MidnightBlue}$Y^h$}
        \put(25,7){\color{MidnightBlue}$\fbox{z}$}
        \put(17.2,79.6){\footnotesize\color{ForestGreen}$-x+ia$}
        \put(20,61.5){\footnotesize$(x,T)$}
        \put(20,41){\footnotesize$(x,0)$}
        \put(27,55){\footnotesize$P$}
        \put(18.5,47){\color{Orchid}\fbox{\footnotesize$(x,s)$}}
        \put(19.0,22.6){\footnotesize\color{MidnightBlue}$x+ia$}
        \put(37.1,91.8){\footnotesize$y'=\tfrac12$}
        \put(37.1,76.8){\footnotesize$y'=a$}
        \put(37.1,66.4){\footnotesize$y'$}
        \put(37.1,57.4){\footnotesize$s=T$}
        \put(37.1,47.9){\footnotesize$s$}
        \put(37.1,3.6){\footnotesize$y=\tfrac12$}
        \put(37.1,38.0){\footnotesize$s=0$}
        \put(37.1,18.9){\footnotesize$y=a$}
        \put(38,31.2){\footnotesize$y$}
        \put(63.3,81.4){\color{Orchid}$C_T$}
        \put(75.5,59.3){\color{ForestGreen}$\overline{Y}^h$}
        \put(60,59.3){\color{ForestGreen}$\fbox{z'}$}
        \put(60.0,65.6){\footnotesize\color{ForestGreen}$-x+ia$}
        \put(60,75){\scriptsize$w'(P)$}
        \put(74,82){\color{Orchid}\fbox{$w'$}}
        \put(87.7,57.0){\footnotesize$\tfrac12$}
        \put(87.7,68.3){\footnotesize$a$}
        \put(87.7,86.8){\footnotesize$a+T$}
        \put(63.3,27.1){\color{Orchid}$C_T$}
        \put(59.5,12.5){\color{MidnightBlue}$Y^h$}
        \put(75,12){\color{MidnightBlue}\fbox{$z$}}
        \put(71.4,17.6){\footnotesize\color{MidnightBlue}$x+ia$}
        \put(71,34.0){\footnotesize$w(P)$}
        \put(58,35){\color{Orchid}\fbox{$w$}}
        \put(87.7,8.6){\footnotesize$\tfrac12$}
        \put(87.7,19.9){\footnotesize$a$}
        \put(87.7,38.9){\footnotesize$a+T$}
        \put(55,51){\small $z'$ continued by $w'$}
        \put(55,3){\small $z$ continued by $w$}
    \end{overpic}
    \caption{Gluing a cylinder $C_T$ to the truncated surfaces $Y^h$ and $\overline{Y}^h$.}\label{Fig-gluing}
\end{figure}

For $T>0$, let
\[
    C_T=\bigl\{(x,s)\in(\R/\Z)\times\R\;\big|\;0\leq s\leq T\bigr\}
\]
be the Euclidean cylinder, with the two holomorphic coordinates
\[
    w,\ w'\colon C_T\longrightarrow\{\zeta\in\C/\Z\mid a\leq\ima(\zeta)\leq a+T\},
\]
\[
    w=x+i(a+s)\qquad\text{and}\qquad w'=-x+i(a+T-s)=-w+i(2a+T).
\]
For each cusp, we glue the circle $\{s=0\}$ of $C_T$ to the horocycle $\{y=a\}$ bounding $Y^h$ with $(x,0)\mapsto x+ia$, in the normalized coordinates $z$ of $Y$ restricted to $Y^h$, and the circle $\{s=T\}$ to the corresponding horocycle $\{y'=a\}$ bounding $\overline Y^h$, by $(x,T)\mapsto -x+ia$, in the normalized coordinates $z'$ of $\overline{Y}^h$. 
Hence, with these gluings, $w$ and $w'$ extend the normalized coordinates $z$ restricted to $Y^h$ and $z'$ of $\overline{Y}^h$ across $\{s=0\}$ and $\{s=T\}$, see Figure~\ref{Fig-gluing}. The gluings are holomorphic, and we denote by $W_T(Y)$ the resulting closed Riemann surface. Note that the cylinder carries the flat metric $|dw|$. 
We call the core curves $\{s=T/2\}$ of the $n$ cylinders the \emph{seams}. Since $Y$ and its mirror $\overline{Y}$ are of type $(g,n)$, the resulting surface $W_T(Y)$ has no puncture and its genus equals $G=2g+(n-1)\geq2$, and its Euler characteristic is $2-2G=2(2-2g-n)$.

\subsubsection{Definition of the marking}
\begin{figure}
    \centering
    \begin{overpic}[width=0.9\linewidth,keepaspectratio]{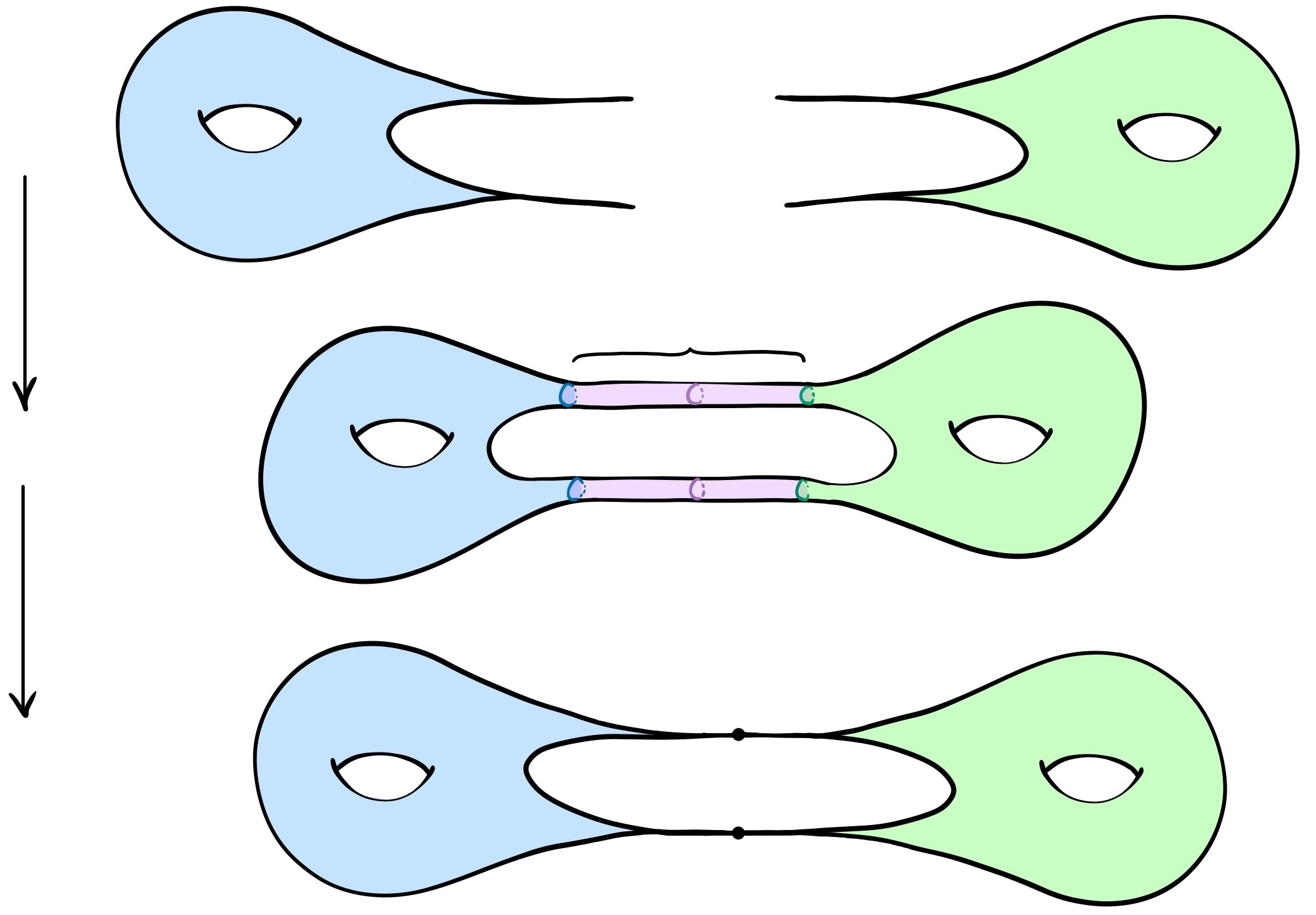}
        \put(15,52){\color{MidnightBlue}$Y$}
        \put(92,52){\color{ForestGreen}$\overline{Y}$}
        \put(-9,47){\small glue $C_T$}
        \put(-8,25){\small $T\to\infty$}
        \put(-2,35){$W_T(Y)$}
        \put(-1,58){$Y,\overline{Y}$}
        \put(23,30){\color{MidnightBlue}$Y^h$}
        \put(79,31){\color{ForestGreen}$\overline{Y}^h$}
        \put(51.5,43.5){$T$}
        \put(51.5,27.5){\color{Orchid}$C_T$}
        \put(48,35){\color{Orchid}seams}
        \put(-2,11){$D(Y)$}
        \put(24,5){\color{MidnightBlue}$Y$}
        \put(87,4){\color{ForestGreen}$\overline{Y}$}
        \put(55,15){$\nu_1$}
        \put(55,4){$\nu_2$}
    \end{overpic}
    \caption{The surfaces $Y$ and $\overline{Y}$ (top), the surface $W_T(Y)$ (middle), and the surface with nodes $D(Y)$ (bottom), here with $g=1$ and $n=2$.}\label{Fig-double}
\end{figure}

Let $K$ be the reference surface $\Sigma$ with an open disc removed around each puncture, and let $\Lambda$ be the closed oriented topological surface obtained by gluing $K$ to a copy $\overline{K}$ of $K$ with the opposite orientation, by an annulus along each pair of corresponding boundary circles. The surface $\Lambda$ is of type $(G,0)$. Let $C\subset \Lambda$ be the multicurve formed by the core curves of the annuli. For a puncture $c$ we denote by $\partial_c K$ the corresponding boundary component of $K$, and we write the annulus at $c$ as $\partial_c K\times[0,1]$, glued to $\partial_c K\subset K$ along $\partial_c K\times\{0\}$ and to the copy of $\partial_c K$ in $\overline{K}$ along $\partial_c K\times\{1\}$. Its \emph{vertical segments} are the segments $\{v\}\times[0,1]$, $v\in\partial_c K$, and its core curve is $\partial_c K\times\{1/2\}$. 
We denote by $d_{\WP}^{\mathrm{aux}}$ the WP distance on $\overline{\teich}(\Lambda)$. Let, $D(Y)\in S(C)$ denotes the marked surface with nodes whose parts are $Y$ and $\overline{Y}$, each cusp of $Y$ being paired with the corresponding cusp of $\overline{Y}$. See Figure~\ref{Fig-double}.

We now define a marking $f_{Y,T} \Lambda \longrightarrow W_T(Y)$. Choose a representative of the marking of $Y$ which maps $K$ onto $Y^h$, and denote by $f_Y\colon K\to Y^h$ its restriction. Near a puncture $c$, the map $f_Y$ sends $\partial_c K$ onto the horocycle $\{y=a\}$, and we define the homeomorphism $\beta_{Y,c}\colon\partial_c K\to\R/\Z$ by
\[
    f_Y(v)=\beta_{Y,c}(v)+ia,\qquad v\in\partial_c K,
\]
in the normalized coordinate $z=x+iy$ at $c$. For $T>0$, let $f_{Y,T}\colon\Lambda\to W_T(Y)$ be equal to $f_Y$ on $K$, to $f_Y$ viewed as a map from $\overline{K}$ onto $\overline{Y}^h$ on $\overline{K}$, and, on the annulus at $c$, in the coordinates $(x,s)$ of the cylinder $C_T$ at $c$,
\[
    f_{Y,T}(v,r)=\bigl(\beta_{Y,c}(v),\,Tr\bigr),\qquad (v,r)\in\partial_cK\times[0,1].
\]
By the gluing rules above, the point $(\beta_{Y,c}(v),0)$ is glued to the point $z=\beta_{Y,c}(v)+ia$ of $Y^h$, which is $f_Y(v)$, and the point $(\beta_{Y,c}(v),T)$ is glued to the point $z'=-\beta_{Y,c}(v)+ia$ of $\overline{Y}^h$, which is $f_Y(v)$ viewed in $\overline{Y}^h$ since $z'=-\bar z$. Hence $f_{Y,T}$ is a well-defined homeomorphism. It restricts to the markings of $Y$ on $K$ and of $\overline{Y}$ on $\overline{K}$, maps vertical segments of the annuli to vertical segments of the cylinders, so that it does not twist along the cylinders, and maps $C$ onto the seams.

\begin{lemma}\label{lem:double-plumbing-limit}
    Let $Y\in\teich_{g,n}$ with $n\geq1$, and let $c=e^{-2\pi a}$. For every $T>0$, the surface $W_T(Y)$, marked by $f_{Y,T}$, is the plumbing of $Y\sqcup\overline{Y}$, the parts of $D(Y)$, at the $n$ pairs of corresponding cusps, with the same parameter
    \[
        t_T=e^{-2\pi(2a+T)}\in(0,c^2)
    \]
    at every pair. Consequently $W_T(Y)\to D(Y)$ in $\overline{\teich}(\Lambda)$ as $T\to\infty$.
\end{lemma}

\begin{proof}
    Fix a pair of corresponding cusps of $Y$ and $\overline{Y}$, with cusp coordinates $q=e^{2\pi iz}$ and $q'=e^{2\pi iz'}$, and set $c=e^{-2\pi a}$. Since $|q|=e^{-2\pi y}$, the annulus $A=\{t_T/c<|q|<c\}$ of the plumbing construction is the annulus $\{a<y<a+T\}$ of $Y$ and the annulus $\{a<y'<a+T\}$ of $\overline{Y}$,

    The coordinates $w$ and $w'$ of the cylinder $C_T$ extend $z$ and $z'$, so the open cylinder is identified with $A$ by $z=w$ and with $A'$ by $z'=w'$. Since $w'=-w+i(2a+T)$, every point of the open cylinder satisfies
    \[
        qq'=e^{2\pi i(w+w')}=e^{-2\pi(2a+T)}=t_T,
    \]
    which is the same rule as in the plumbing construction. Hence $W_T(Y)=Y^h\cup C_T\cup\overline{Y}^h$ is obtained from $Y\sqcup\overline{Y}$ by removing the discs $\{|q|\leq t_T/c\}$ and $\{|q'|\leq t_T/c\}$ and identifying $A$ with $A'$ by $qq'=t_T$, that is, it is the plumbing with parameter $t_T$ at every pair of cusps. Moreover its marking $f_{Y,T}$ restricts to the markings of $Y$ and $\overline{Y}$ on $K_c=Y^h\sqcup\overline{Y}^h$.

    Since $a=4$, we have $c=e^{-8\pi}\leq e^{-\pi}$ and $t_T=c^2e^{-2\pi T}\in(0,c^2)$. As $T\to\infty$, $t_T\to0$, and Proposition~\ref{prop:plumbing} gives $W_T(Y)\to D(Y)$ in $\overline{\teich}(\Lambda)$.
\end{proof}

The next lemma bounds the WP length along the path $T\mapsto W_T(Y)$.

\begin{lemma}\label{lem:opening}
    For every $T>0$,
    \[
        \left\|\frac{d}{dT}W_T(Y)\right\|_{\WP}^2\leq\frac{2\pi n}{T^3}
        \qquad\text{and}\qquad
        d_{\WP}^{\mathrm{aux}}\bigl(W_T(Y),D(Y)\bigr)\leq2\sqrt{\frac{2\pi n}{T}}.
    \]
\end{lemma}

\begin{proof}
    Fix $T>0$. Denote by $C_{T,1},\dots,C_{T,n}$ the cylinders of $W_T(Y)$, one for each cusp of $Y$, and by $\zeta_{T,k}=x+is$, with $x\in\R/\Z$ and $0\leq s\leq T$, the coordinate of $C_{T,k}$. 
    On $C_{T,k}$ the coordinates $w=\zeta_{T,k}+ia$ and $w'=-\zeta_{T,k}+i(a+T)$ extend the restriction of the normalized coordinates $z$ of $Y$ to $Y^h$ and $z'$ of $\overline{Y}^h$ across the gluing circles $\{s=0\}$ and $\{s=T\}$. Let $v_T\colon[0,T]\to[0,1]$ be the piecewise linear function
    \[
        v_T(s)=
        \begin{cases}
            0, & 0\leq s\leq T/4,\\
            2s/T-1/2, & T/4\leq s\leq 3T/4,\\
            1, & 3T/4\leq s\leq T,
        \end{cases}
    \]
    and for $|\varepsilon|<T/2$ let $\phi_\varepsilon(s)=s+\varepsilon v_T(s)$. 
    Define $F_\varepsilon\colon W_T(Y)\to W_{T+\varepsilon}(Y)$ by
    \[
        \zeta_{T+\varepsilon,k}\bigl(F_\varepsilon(x+is)\bigr)=x+i\phi_\varepsilon(s)\qquad\text{on }C_{T,k},\ k=1,\dots,n,
    \]
    and by the identity on $Y^h$ and $\overline{Y}^h$, see Figure~\ref{Fig-opening}. On $\{0\leq s\leq T/4\}$ we have $\phi_\varepsilon(s)=s$, hence $w\circ F_\varepsilon=w$. On $\{3T/4\leq s\leq T\}$ we have $\phi_\varepsilon(s)=s+\varepsilon$, and since $w'=-\zeta_{T+\varepsilon,k}+i(a+T+\varepsilon)$ on $C_{T+\varepsilon,k}$, we obtain:
    \[
        w'\bigl(F_\varepsilon(x+is)\bigr)=-x+i\bigl(a+T+\varepsilon-(s+\varepsilon)\bigr)=-x+i(a+T-s)=w'(x+is).
    \]
    \begin{figure}[h]
    \centering
    \begin{overpic}[width=0.9\linewidth,keepaspectratio]{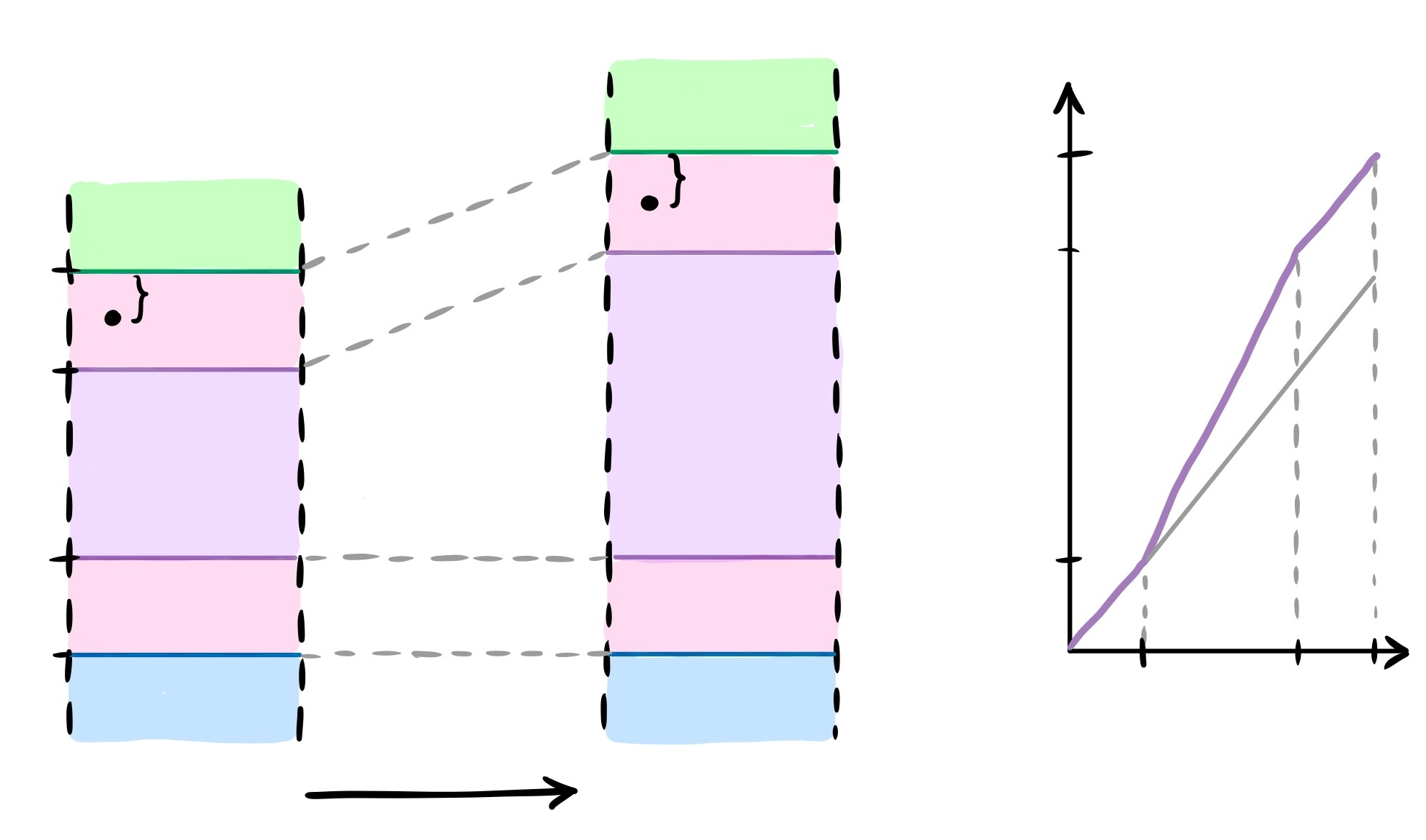}
        \put(1,12.3){$0$}
        \put(0,18.8){$\tfrac{T}{4}$}
        \put(-1,32.2){$\tfrac{3T}{4}$}
        \put(12,37){\footnotesize$T-s$}
        \put(49,45.5){\footnotesize$T-s$}
        \put(2,25){\color{Orchid}$s$}
        \put(0,39.2){$T$}
        \put(9,26){\color{Orchid}$M_{T,k}$}
        \put(45,29){\color{Orchid}$F_\varepsilon(M_{T,k})$}
        \put(11,3){$C_{T,k}$}
        \put(43,3){$C_{T+\varepsilon,k}$}
        \put(28,5.5){$F_\varepsilon$}
        \put(73,11){$0$}
        \put(71,18.8){$\tfrac{T}{4}$}
        \put(64,40.3){$\tfrac{3T}{4}+\varepsilon$}
        \put(66,47.2){$T+\varepsilon$}
        \put(79,8.3){$\tfrac{T}{4}$}
        \put(89.2,8.3){$\tfrac{3T}{4}$}
        \put(95,8.3){$T$}
        \put(98.5,10.5){\footnotesize$s$}
        \put(76.5,51){\footnotesize\color{Orchid}$\phi_\varepsilon(s)$}
    \end{overpic}
    \caption{The map $F_\varepsilon\colon W_{T}\to W_{T+\varepsilon}$ near one cylinder.}\label{Fig-opening}
    \end{figure}

    Thus, near each gluing circle, $F_\varepsilon$ is the identity in a coordinate extending the normalized coordinate of $Y^h$ and $\overline{Y}^h$. Hence $F_\varepsilon$ is a well-defined homeomorphism compatible with the gluings. 
    Moreover, $F_\varepsilon$ is compatible with the markings, that is, $F_\varepsilon\circ f_{Y,T}$ is isotopic to $f_{Y,T+\varepsilon}$. Indeed, both maps are equal to $f_Y$ on $K$ and $\overline{K}$, since $F_\varepsilon$ is the identity on $Y^h$ and $\overline{Y}^h$. On the annulus at $c$, both send $(v,r)$ to the vertical segment $\{x=\beta_{Y,c}(v)\}$ of $C_{T+\varepsilon}$, at the heights $s=\phi_\varepsilon(Tr)$ and $s=(T+\varepsilon)r$ respectively.
    Both heights increase from $0$ to $T+\varepsilon$ as $r$ goes from $0$ to $1$. Hence the two maps agree on the boundary circles of the annulus and differ only by an increasing reparametrization of each vertical segment fixing its endpoints, so they are isotopic.
    Hence $F_\varepsilon$ maps the marked surface $W_T(Y)$ to the marked surface $W_{T+\varepsilon}(Y)$, and its Beltrami coefficient represents the variation at $W_T(Y)$ of the path $\tau\mapsto W_\tau(Y)$ in $\teich(\Lambda)$ for $\tau>0$.
    
    We compute the Beltrami coefficient of $F_\varepsilon$. On $C_{T,k}$ write $\zeta=\zeta_{T,k}=x+is$, so that $F_\varepsilon(\zeta)=x+i\phi_\varepsilon(s)$. Since $\partial_\zeta=\frac12(\partial_x-i\partial_s)$ and $\partial_{\bar\zeta}=\frac12(\partial_x+i\partial_s)$, we have $\partial_\zeta F_\varepsilon=\frac12\bigl(1+\phi_\varepsilon'(s)\bigr)$ and $\partial_{\bar\zeta}F_\varepsilon=\frac12\bigl(1-\phi_\varepsilon'(s)\bigr)$, hence
    \[
        \frac{\partial_{\bar\zeta}F_\varepsilon}{\partial_\zeta F_\varepsilon}=\frac{1-\phi_\varepsilon'(s)}{1+\phi_\varepsilon'(s)}.
    \]
    Let $M_{T,k}=\{T/4\leq s\leq 3T/4\}\subset C_{T,k}$ be the middle half of $C_{T,k}$ and $M_T=M_{T,1}\cup\dots\cup M_{T,n}$. Since $\phi_\varepsilon'=1+2\varepsilon/T$ on $M_T$ and $\phi_\varepsilon'=1$ on the rest of the cylinders, the map $F_\varepsilon$ has Beltrami coefficient
    \begin{equation}\label{eq-beltrami-coefficient}
        \mu_\varepsilon=-\frac{\varepsilon}{T+\varepsilon}\,\frac{\overline{d\zeta_{T,k}}}{d\zeta_{T,k}}\quad\text{on }M_{T,k},\ \text{for } k=1,\dots,n,
    \end{equation}
    and $\mu_\varepsilon=0$ outside $M_T$. Since $F_\varepsilon$ is piecewise affine and $|\mu_\varepsilon|<1$, it is quasiconformal. The value $\varepsilon=0$ corresponds to $W_T(Y)$, since $F_0$ is the identity. Taking the derivative of $-\varepsilon/(T+\varepsilon)$ with respect to $\varepsilon$ and setting $\varepsilon=0$ gives $-1/T$. Thus the tangent vector at $W_T(Y)$ to the path $\tau\mapsto W_\tau(Y)$ is represented by the Beltrami differential
    \[
        \mu_F(T)=-\frac{1}{T}\,\frac{\overline{d\zeta_T}}{d\zeta_T},
    \]
    where $\frac{\overline{d\zeta_T}}{d\zeta_T}=\frac{\overline{d\zeta_{T,k}}}{d\zeta_{T,k}}$ on $M_{T,k}$, $k=1,\dots,n$, and zero outside $M_T$. By~\eqref{eq:beltrami-bound},
    \begin{equation}\label{eq-WP-norm-beltrami}
         \left\|\frac{d}{dT}W_T(Y)\right\|_{\WP}^2\leq\int_{W_T(Y)}|\mu_F(T)|^2\,dA=\frac{1}{T^2}\sum_{k=1}^n\area(M_{T,k}),
    \end{equation}
    where $dA$ and $\area$ refer to the hyperbolic metric of $W_T(Y)$. 
    
    To bound $\area(M_{T,k})$, we compare the hyperbolic metric of $W_T(Y)$ with the complete hyperbolic metric of the open cylinder $C_{T,k}^\circ=\{0<s<T\}\subset C_{T,k}$, as in~\cite[Section 2]{Wolpert1975}. The open cylinder $C_{T,k}^\circ$ is the quotient of the strip $\{\zeta\in\C\mid 0<\ima(\zeta)<T\}$ by the translation $\zeta\mapsto\zeta+1$, and since the strip is simply connected, it is the universal cover of $C_{T,k}^\circ$. The map $\zeta \mapsto \omega=e^{\pi\zeta/T}=e^{\pi x/T}e^{i\pi s/T}$ sends this strip biholomorphically onto $\Hh$ because $\arg\omega=\pi s/T\in(0,\pi)$. 
    Since $d\omega=\frac{\pi}{T}\,\omega\,d\zeta$ and $\ima(\omega)=e^{\pi x/T}\sin(\pi s/T)=|\omega|\sin(\pi s/T)$, the pull-back of the hyperbolic metric $|d\omega|/\ima(\omega)$ of $\Hh$ is
    \[
        \frac{|d\omega|}{\ima(\omega)}=\frac{\frac{\pi}{T}\,|\omega|\,|d\zeta|}{|\omega|\sin(\pi s/T)}=\lambda_T(s)\,|d\zeta|, \quad
        \text{where }
        \lambda_T(s)=\frac{\pi}{T}\csc\Bigl(\frac{\pi s}{T}\Bigr),\quad 0<s<T.
    \]
    Since $\lambda_T$ only depends only on $s$ (recall that $T$ is fixed here), it is invariant under $\zeta\mapsto\zeta+1$ and descends to the complete hyperbolic metric of $C_{T,k}^\circ$.
    In the coordinate $\zeta=x+is$, write the hyperbolic metric of $W_T(Y)$ as $\rho(\zeta)^2|d\zeta|^2$, so that its area form is $dA=\rho(\zeta)^2\,dx\,ds$.
    The inclusion $C_{T,k}^\circ\hookrightarrow W_T(Y)$ is holomorphic, so it lifts to a holomorphic map from the strip to the universal cover $\Hh$ of $W_T(Y)$. By the Schwarz--Pick lemma this lift does not increase hyperbolic lengths, and since covering maps are local isometries we obtain
    \[
        \rho(\zeta)\leq\lambda_T(s)\qquad\text{on }C_{T,k}^\circ.
    \]
    Hence
    \begin{equation}\label{eq-bound-area-cylinder}
        \area(M_{T,k})\leq\int_0^1\int_{T/4}^{3T/4}\frac{\pi^2}{T^2}\csc^2\left(\frac{\pi s}{T}\right)ds\,dx=\frac{\pi}{T}\Bigl[-\cot\Bigl(\frac{\pi s}{T}\Bigr)\Bigr]_{T/4}^{3T/4}=\frac{2\pi}{T},
    \end{equation}
    and by combining~\eqref{eq-WP-norm-beltrami} and~\eqref{eq-bound-area-cylinder} we obtain:
    \[
        \left\|\frac{d}{dT}W_T(Y)\right\|_{\WP}\leq \left(\frac{1}{T^2}\sum_{k=1}^n\area(M_{T,k}) \right)^{1/2} \leq \sqrt{\frac{2\pi n}{T^3}}.
    \]
    Consequently, for $T\leq T'$,
    \begin{align*}
        d_{\WP}^{\mathrm{aux}}\bigl(W_T(Y),W_{T'}(Y)\bigr) &\leq \int_T^{T'}\left\|\frac{d}{dt}W_t(Y)\right\|_{\WP} dt \leq\int_T^{T'}\sqrt{2\pi n}\,t^{-3/2}\,dt\\
        &=2\sqrt{\frac{2\pi n}{T}} - 2\sqrt{\frac{2\pi n}{T'}}\leq 2\sqrt{\frac{2\pi n}{T}}.
    \end{align*}
    Hence $(W_t(Y))_{t\geq T}$ is a Cauchy family for the WP metric. By Lemma~\ref{lem:double-plumbing-limit} its limit is $D(Y)$. Letting $T'\to\infty$ gives the distance estimate
    \[
        d_{\WP}^{\mathrm{aux}}\bigl(W_T(Y),D(Y)\bigr)\leq2\sqrt{\frac{2\pi n}{T}}.
    \]
\end{proof}

\subsection{Lengths of the doubled arcs}\label{section-def-de}
Set
\[
    T_0=\frac{\pi a}{2}=2\pi,\qquad b=a+\frac{T_0}{2}=a+\frac{\pi a}{4}=4+\pi,
\]
and write $W_{T_0}(Y)$. On each cylinder $C_{T_0,k}$ of $W_{T_0}(Y)$ we use the coordinate $w=x+iy$ of Section~\ref{subsec:opening}, with $y=a+s\in[a,a+T_0]$. It extends across the circle $\{y=a\}$ the restriction of the normalized coordinate $z=x+iy$ of $Y$ to the band $\{1/2<y\leq a\}\subset Y^h$. 
The seam $\{s=T_0/2\}$ of $C_{T_0,k}$ is the circle $\{y=b\}$, and the mirror coordinate is $w'=-x+iy'$ with $y'=a+T_0-s=2b-y$.

\begin{figure}[h]
    \centering
    \begin{overpic}[width=0.9\linewidth,keepaspectratio]{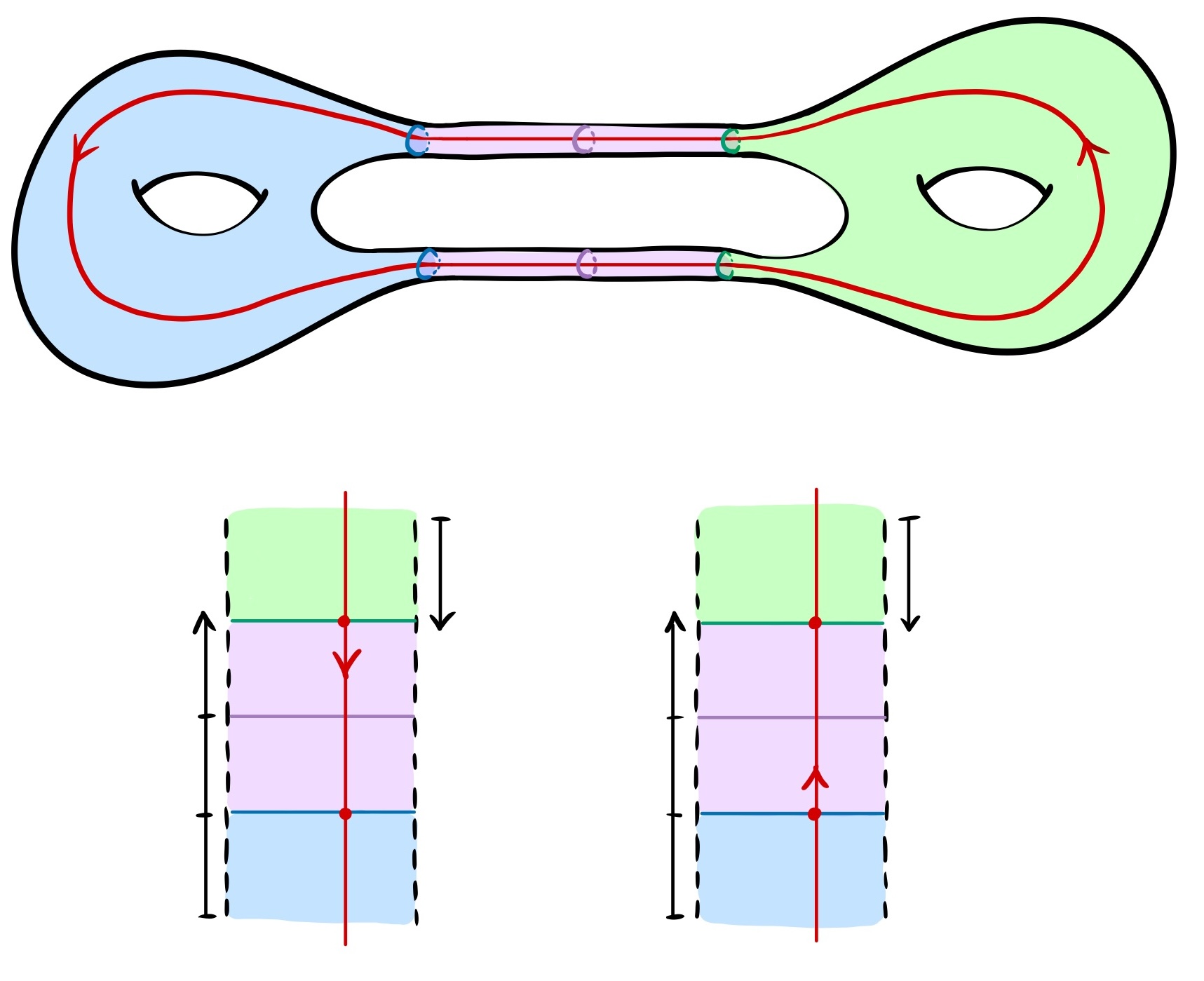}
        \put(6,52){\color{MidnightBlue}$Y^h$}
        \put(91,72){\color{ForestGreen}$\overline{Y}^h$}
        \put(85,61){\color{red}$d(e)$}
        \put(39,73.5){\small cylinder at $c_1$}
        \put(39,53.5){\small cylinder at $c_2$}
        \put(30.5,8){\footnotesize\color{MidnightBlue}$z_1$}
        \put(30.5,17){\footnotesize$w_1$}
        \put(30.5,25){\footnotesize$w'_1$}
        \put(30.5,34){\footnotesize\color{ForestGreen}$z'_1$}
        \put(22,12){\tiny\color{red}$x_1{+}ia$}
        \put(20,30.5){\tiny\color{red}$-x_1{+}ia$}
        \put(14,4.8){\footnotesize$\tfrac12$}
        \put(14.5,13.3){\footnotesize$a$}
        \put(14.5,21.3){\footnotesize$b$}
        \put(8,29.3){\footnotesize$a+T_0$}
        \put(16.2,33.2){\footnotesize$y$}
        \put(38.3,37.3){\footnotesize$\tfrac12$}
        \put(38.3,28.8){\footnotesize$a$}
        \put(37,33.5){\footnotesize$y'$}
        \put(69.5,8){\footnotesize\color{MidnightBlue}$z_2$}
        \put(69.5,17){\footnotesize$w_2$}
        \put(69.5,25){\footnotesize$w'_2$}
        \put(69.5,34){\footnotesize\color{ForestGreen}$z'_2$}
        \put(61,12){\tiny\color{red}$x_2{+}ia$}
        \put(59.0,30.5){\tiny\color{red}$-x_2{+}ia$}
        \put(53,4.8){\footnotesize$\tfrac12$}
        \put(53.5,13.3){\footnotesize$a$}
        \put(53.5,21.3){\footnotesize$b$}
        \put(47,29.3){\footnotesize$a+T_0$}
        \put(55.2,33.2){\footnotesize$y$}
        \put(77.3,37.3){\footnotesize$\tfrac12$}
        \put(77.3,28.8){\footnotesize$a$}
        \put(76,33.5){\footnotesize$y'$}
        \put(19,-1){\small near the cusp $c_1$}
        \put(58,-1){\small near the cusp $c_2$}
    \end{overpic}\vspace*{0.3cm}
    \caption{The doubled arc $d(e)$ on $W_{T_0}(Y)$, for an arc $e$ of $\T$ from the cusp $c_1$ to the cusp $c_2$.}\label{Fig-darc}
\end{figure}

Let $e$ be an arc of $\T$, realized as a geodesic on $Y$, and let $c_1$ and $c_2$ be the cusps at which $e$ ends, with $c_1=c_2$ allowed. For $j=1,2$, denote by $z_j=x+iy$ the normalized coordinate of $Y$ at $c_j$, by $z'_j=-\bar z_j$ the normalized coordinate of $\overline{Y}$ at $c_j$, and by $w_j=x+iy$ and $w'_j=-x+iy'$ the coordinates of the cylinder of $W_{T_0}(Y)$ glued at $c_j$. Then, $e\cap Y^h$ is a geodesic segment of length $\delta_e(Y)$, orthogonal to $\partial Y^h$, whose endpoints are $z_1=x_1+ia$ on the horocycle $\{y=a\}$ at $c_1$ and $z_2=x_2+ia$ on the horocycle $\{y=a\}$ at $c_2$. Since $e$ meets the standard cusp neighborhoods only along terminal rays~\cite[proof of Lemma 8.8]{abadie2026}, near its endpoint at $c_j$ the segment $e\cap Y^h$ is the vertical segment $\{x=x_j,\ 1/2<y\leq a\}$ in the coordinate $z_j$. We denote by $d(e)$ the closed curve on $W_{T_0}(Y)$ obtained by following
\begin{enumerate}
    \item the segment $e\cap Y^h$, from $z_1=x_1+ia$ to $z_2=x_2+ia$;
    \item the vertical segment $\{x=x_2,\ a\leq y\leq a+T_0\}$ of the cylinder at $c_2$, in the coordinate $w_2$, from $w_2=x_2+ia$ to $w'_2=-x_2+ia$;
    \item the mirror image of $e\cap Y^h$ in $\overline{Y}^h$, from $z'_2=-x_2+ia$ to $z'_1=-x_1+ia$;
    \item the vertical segment $\{x=x_1,\ a\leq y\leq a+T_0\}$ of the cylinder at $c_1$, in the coordinate $w_1$, from $w'_1=-x_1+ia$ to $w_1=x_1+ia$.
\end{enumerate}
See Figure~\ref{Fig-darc}. If $c_1=c_2$, the two vertical segments lie in the same cylinder, at the distinct positions $x_1\neq x_2$.

\begin{lemma}[Doubled arc estimate]\label{lem:arcs}
    Let $Y\in\teich_{g,n}$ with $n\geq1$, let $\T$ be an ideal triangulation of $Y$, and let $T_0=2\pi$. For every arc $e$ of $\T$, the length of the closed geodesic of $W_{T_0}(Y)$ freely homotopic to $d(e)$ satisfies
    \[
        \ell_{W_{T_0}(Y)}(d(e))\leq2\sqrt2\,\delta_e(Y)+4\log(1+\sqrt2).
    \]
\end{lemma}

\begin{proof}
    We construct a conformal metric $g_*$ on $W_{T_0}(Y)$ which is larger than the hyperbolic metric $g_W$ of $W_{T_0}(Y)$, and we bound the $g_*$-length of $d(e)$. In a local coordinate $\zeta$ we write $g_W=\rho(\zeta)^2|d\zeta|^2$ and $g_*=\rho_*(\zeta)^2|d\zeta|^2$. Let $J\subset W_{T_0}(Y)$ be the union of the $2n$ gluing circles $\{y=a\}$ and $\{y'=a\}$, that is $J=\partial Y^h\cup\partial\overline{Y}^h$.

    On $Y^h$ and $\overline{Y}^h$ we set $g_*=2g_Y$ and $g_*=2g_{\overline{Y}}$, where $g_Y$ and $g_{\overline{Y}}$ are the hyperbolic metrics of $Y$ and $\overline{Y}$. On each cylinder we set $g_*=\rho_*(y)^2|dw|^2$, where
    \[
        \rho_*(y)=\frac{1}{a}\sec\left(\frac{b-y}{a}\right),\qquad a\leq y\leq a+T_0.
    \]
    Note that for $a\leq y\leq a+T_0$ we have $|b-y|/a\leq T_0/(2a)=\pi/4<\pi/2$ thus $\cos\bigl(\frac{b-y}{a}\bigr)\geq\cos(\pi/4)>0$ and so the density $\rho_*$ is well defined and positive. 
    Since $2b=2a+T_0$, for every $y\in[a,a+T_0]$ we also have $2b-y\in[a,a+T_0]$, indeed the map $y\mapsto2b-y$ is the reflection of the cylinder across the seam $\{y=b\}$, and it exchanges the two gluing circles $\{y=a\}$ and $\{y=a+T_0\}=\{y'=a\}$. Since $\sec$ is even, 
       \[
            \rho_*(2b-y)=\frac1a\sec\left(\frac{b-(2b-y)}{a}\right)=\frac1a\sec\left(\frac{y-b}{a}\right)=\frac1a\sec\left(\frac{b-y}{a}\right)=\rho_*(y).
        \]
    Moreover the mirror coordinate is $w'=-x+iy'$ with $y'=2b-y$ hence $dw'=-dw$ and $|dw'|=|dw|$. Therefore
    \[
        g_*=\rho_*(y)^2|dw|^2=\rho_*(2b-y')^2|dw'|^2=\rho_*(y')^2|dw'|^2.
    \]
    Thus, by construction, $g_*$ is given by the same formula in the mirror coordinate. On the band $\{1/2<y\leq a\}\subset Y^h$, the hyperbolic metric of $Y$ in the normalized coordinate is $g_Y=|dz|^2/y^2$, hence
    \[
        g_*=2g_Y=\frac{2|dz|^2}{y^2}=\Bigl(\frac{\sqrt2}{y}\Bigr)^2|dz|^2.
    \]
    Since $w=z$ on this band, in the coordinate $w$ extended by $z$ the metric $g_*$ is $\rho_*(y)^2|dw|^2$ on $\{1/2<y\leq a+T_0\}$, with
    \[
        \rho_*(y)=
        \begin{cases}
            \sqrt2/y, & 1/2<y\leq a,\\[2pt]
            \frac1a\sec\bigl(\frac{b-y}{a}\bigr), & a\leq y\leq a+T_0,
        \end{cases}
    \]
    and, by symmetry, $g_*=\rho_*(y')^2|dw'|^2$ near $\{y'=a\}$, with $\rho_*(y')=\sqrt2/y'$ on the band $\{1/2<y'\leq a\}\subset\overline{Y}^h$. At $y=a$ both expressions give the value $\sqrt2/a$, so $\rho_*$ is continuous across $\{y=a\}$, and by the symmetry above also across $\{y'=a\}$.

    Hence $g_*$ is a well-defined continuous metric on $W_{T_0}(Y)$. The surface $W_{T_0}(Y)$ is the union of the compact surfaces with boundary $Y^h$, $\overline{Y}^h$ and $C_{T_0,1},\dots,C_{T_0,n}$, glued along $J$, and we need the following regularity of $\rho_*$ at $J$ in the next part of the proof.

    \begin{claim}\label{claim-rho-J}
        On each of the surfaces $Y^h$, $\overline{Y}^h$ and $C_{T_0,k}$, the function $\log\rho_*$ is the restriction of a smooth function defined on a neighborhood of this surface in $W_{T_0}(Y)$. Moreover $\log\rho_*$ is of class $C^1$ on $W_{T_0}(Y)$, that is, across each circle of $J$ the two expressions of $\log\rho_*$ have the same value and the same first derivatives.
    \end{claim}

    \begin{proof}[Proof of Claim~\ref{claim-rho-J}]
        By symmetry it suffices to consider $Y^h$, the cylinders, and the circles $\{y=a\}$. On $Y^h$, $\rho_*$ is the density of $2g_Y$, and near $\partial Y^h$ it is $\sqrt2/y$, which is smooth for $y>0$. On $C_{T_0,k}$, the expression $\frac1a\sec\bigl(\frac{b-y}{a}\bigr)$ is smooth and positive for $|b-y|<\pi a/2$, an interval containing $[a,a+T_0]$ since $b-a=\pi a/4$. This proves the first statement.
        At $y=a$ the exressions for $\rho_*$ on both sides agree and take the value $\sqrt2/a$. 
        Both expressions depend only on $y$, so $\partial_x\log\rho_*=0$ on both sides, and it remains to compare the $y$-derivatives at $y=a$. On the side of $Y^h$, $(\log\rho_*)'(y)=-1/y$, which equals $-1/a$. On the side of the cylinder, $(\log\rho_*)'(y)=-\frac1a\tan\bigl(\frac{b-y}{a}\bigr)$, which equals $-\frac1a\tan(\pi/4)=-1/a$.
    \end{proof}
    
    The curvature of $g_*=\rho_*^2|d\zeta|^2$ is $K_*=-\rho_*^{-2}\Delta\log\rho_*$, where $\Delta=\partial_x^2+\partial_y^2$ is the Euclidean Laplacian in the coordinate $\zeta=x+iy$. On $Y^h$ and $\overline{Y}^h$ we have $K_*=-1/2$, since $g_*$ is twice a hyperbolic metric. On the cylinders $\rho_*$ only depends on $y$, and
    \[
        \Delta\log\rho_*=(\log\rho_*)''(y)=\frac{1}{a^2}\sec^2\left(\frac{b-y}{a}\right)=\rho_*(y)^2,
    \]
    so $K_*=-1$. Hence $K_*\geq-1$ on $W_{T_0}(Y)\smallsetminus J$, that is
    \begin{equation}\label{eq:comparison-curvature}
        \Delta\log\rho_*\leq\rho_*^2\qquad\text{on }W_{T_0}(Y)\smallsetminus J.
    \end{equation}
    Note that on $J$ the function $\log\rho_*$ is only of class $C^1$: at $y=a$ its second derivative is $1/a^2$ on the side of $Y^h$ and $\frac{1}{a^2}\sec^2(\pi/4)=2/a^2$ on the side of the cylinder, so $\Delta\log\rho_*$ is not defined on $J$. However~\eqref{eq:comparison-curvature} holds on each component of $W_{T_0}(Y)\smallsetminus J$.
     
    The hyperbolic metric of $W_{T_0}(Y)$ has curvature $-1$, that is $\Delta\log\rho=\rho^2$ in every local coordinate. Both $g_W$ and $g_*$ are conformal metrics on the Riemann surface $W_{T_0}(Y)$, and so we write \[g_W=e^{2u}g_*,\] that is $\rho^2=\rho_*^2 e^{2u}$. 
    In a holomorphic local coordinate $\zeta$, write $g_W=\rho(\zeta)^2|d\zeta|^2$ and $g_*=\rho_*(\zeta)^2|d\zeta|^2$. Changing to another holomorphic coordinate $\xi$ multiplies both densities by $|d\zeta/d\xi|$, so the ratio $\rho/\rho_*$ does not depend on the coordinate, and so $u=\log(\rho/\rho_*)$ is a well-defined function on $W_{T_0}(Y)$.
    Since $\log\rho$ is smooth on $W_{T_0}$ (in local coordinates), Claim~\ref{claim-rho-J} shows that $u=\log\rho-\log\rho_*$ is of class $C^1$ on $W_{T_0}$, and that its restriction to each of the surfaces $Y^h$, $\overline{Y}^h$ and $C_{T_0,1},\dots,C_{T_0,n}$ extends to a smooth function on a neighborhood of this surface. On $J$, $u$ is only $C^1$. 
    \begin{claim}\label{claim-u-negative}
        We have $u\leq0$ on $W$, that is, $g_W\leq g_*$.
    \end{claim}
    \begin{proof}[Proof of Claim~\ref{claim-u-negative}]
        By~\eqref{eq:comparison-curvature}, on $W_{T_0}(Y)\smallsetminus J$,
        \begin{equation}\label{eq:comparison-weak}
            \Delta u=\rho^2-\Delta\log\rho_*\geq\rho^2-\rho_*^2=\rho_*^2\bigl(e^{2u}-1\bigr).
        \end{equation}
        Since $W_{T_0}(Y)$ is compact, $u$ attains its maximum $M$ at some point $p\in W_{T_0}(Y)$. Assume by contradiction that $M>0$.
        \begin{itemize}
            \item If $p\notin J$, then $u$ is smooth near $p$ and $\Delta u(p)\leq0$, since $p$ is an interior maximum. On the other hand~\eqref{eq:comparison-weak} gives $\Delta u(p)\geq\rho_*(p)^2(e^{2M}-1)>0$, a contradiction.
            \item If $p\in J$, by symmetry we may assume that $p$ lies on a circle $\{y=a\}$, say $p=x_0+ia$ in the coordinate $w$ extended by $z$. Let $u_+$ be the restriction of $u$ to the cylinder $C_{T_0,k}$ containing $p$, which is smooth on $C_{T_0,k}$, including its boundary circle $\{y=a\}$. Since $u$ is $C^1$ and has a maximum at $p$, we have $\partial_x u(p)=\partial_y u(p)=0$. The function $x\mapsto u_+(x,a)$ has a maximum at $x_0$, so $\partial_x^2u_+(p)\leq0$. The function $y\mapsto u_+(x_0,y)$, $y\geq a$, has a maximum and a vanishing derivative at $y=a$, so $\partial_y^2u_+(p)\leq0$, since otherwise $u_+(x_0,y)>M$ for $y>a$ close to $a$. Hence $\Delta u_+(p)\leq0$. On the other hand, letting $y\to a^+$ in~\eqref{eq:comparison-weak} gives $\Delta u_+(p)\geq\rho_*(p)^2(e^{2M}-1)>0$, a contradiction.
        \end{itemize}
        Therefore $u\leq0$ and so $g_W\leq g_*$.
    \end{proof}

    The closed geodesic freely homotopic to $d(e)$ minimizes the $g_W$-length in its free homotopy class, hence by the claim above
    \[
        \ell_{W_{T_0}(Y)}(d(e))\leq\ell_{g_W}(d(e))\leq\ell_{g_*}(d(e)).
    \]
    Since $g_*=2g_Y$ on $Y^h$ and $g_*=2g_{\overline{Y}}$ on $\overline{Y}^h$, the segment $e\cap Y^h$ and its mirror image each have $g_*$-length $\sqrt2\,\delta_e(Y)$. For $j=1,2$, the vertical segment $\{x=x_j,\ a\leq y\leq a+T_0\}$ of the cylinder at $c_j$, where $y=\ima(w_j)$, has $g_*$-length,
    \[
        \int_a^{a+T_0}\frac{1}{a}\sec\left(\frac{b-y}{a}\right)dy=\int_{-\pi/4}^{\pi/4}\sec t\,dt=2\log(1+\sqrt2),
    \]
    Summing these four contributions gives $\ell_{g_*}(d(e))=2\sqrt2\,\delta_e(Y)+4\log(1+\sqrt2)$, which proves the lemma.
\end{proof}

\subsection{Compatible markings and a common maximally noded surface}\label{subsec:markings}
In this section we construct a pants decomposition $\mathcal{P}_{\T}$ of the reference surface $\Lambda$ whose curves are mapped by the marking $f_{Y,T_0}$ to curves freely homotopic to the doubled arcs $d(e)$ on $W_{T_0}$ for $e\in\T$. Pinching $\mathcal{P}_{\T}$ gives a maximally noded surface $P_{\T}$ in $\overline{\teich}(\Lambda)$. Recall from Section~\ref{subsec:opening} that $\Lambda$ is obtained from $K$ and its mirror copy $\overline{K}$ by gluing the annuli $\partial_cK\times[0,1]$, and that $W_T(Y)$ is marked by the homeomorphism $f_{Y,T}\colon\Lambda\to W_T(Y)$, which maps $K$ onto $Y^h$, $\overline{K}$ onto $\overline{Y}^h$, and each vertical segment of an annulus onto a vertical segment of the corresponding cylinder. We first construct on $\Lambda$ the curves corresponding to the doubled arcs.

Fix pairwise disjoint simple arcs $a_e\subset K$, $e\in\T$, with pairwise distinct endpoints on $\partial K$, such that each $a_e$ \emph{represents} $e$, that is, extending $a_e$ at each of its endpoints by a ray of the removed punctured disc to the corresponding puncture gives an arc of $\Sigma$ in the isotopy class of $e$. For $e\in\T$, let $d_\Lambda(e)\subset\Lambda$ be the closed curve formed by $a_e$, its copy $\bar a_e\subset\overline{K}$, and the two vertical segments of the annuli through the endpoints of $a_e$. The curves $d_\Lambda(e)$ are well defined up to isotopy. 

\begin{lemma}\label{lem:markings}
    Let $\T$ be an ideal triangulation of $\Sigma$, and let $d_\Lambda(e)\subset\Lambda$, $e\in\T$, be the curves defined above.
    \begin{enumerate}
        \item For every $Y\in\teich(\Sigma)$ and every $e\in\T$, the curve $f_{Y,T_0}(d_\Lambda(e))$ is freely homotopic in $W_{T_0}(Y)$ to the curve $d(e)$ of Section~\ref{section-def-de}.
        \item The curves $d_\Lambda(e)$, $e\in\T$, form a pants decomposition $\mathcal P_{\T}$ of $\Lambda$. In particular, the stratum $S(\mathcal P_{\T})\subset\overline{\teich}(\Lambda)$ is a single maximally noded surface $P_{\T}$.
    \end{enumerate}
\end{lemma}

\begin{proof}
    \begin{enumerate}
        \item Let $c_1$ and $c_2$ be the cusps at which $e$ ends, as in the definition of $d(e)$ (see the beginning of Section~\ref{section-def-de}). We orient $e$ from $c_1$ to $c_2$, and $a_e$ accordingly, from its endpoint $v_1\in\partial_{c_1}K$ to its endpoint $v_2\in\partial_{c_2}K$. By construction, $f_{Y,T_0}(d_\Lambda(e))$ follows the arc $f_Y(a_e)\subset Y^h$, the vertical segment of the cylinder at $f_Y(v_2)$, the mirror image of $f_Y(a_e)$ in $\overline{Y}^h$, and the vertical segment of the cylinder at $f_Y(v_1)$, while $d(e)$ is obtained in the same way from the segment $e\cap Y^h$. Since $a_e$ represents $e$, there is a homotopy $(\gamma_t)_{t\in[0,1]}$ from $\gamma_0=f_Y(a_e)$ to $\gamma_1=e\cap Y^h$ through paths in $Y^h$ from the horocycle $\{y=a\}$ at $c_1$ to the horocycle $\{y=a\}$ at $c_2$. With the coordinates $z_j$, $z'_j$, $w_j$ and $w'_j$, $j=1,2$, used in the definition of $d(e)$, we write the endpoints of $\gamma_t$ as
        \[
            z_1=x_1(t)+ia\qquad\text{and}\qquad z_2=x_2(t)+ia,
        \]
        where $x_1,x_2\colon[0,1]\to\R/\Z$ are continuous. For $t\in[0,1]$, let $\delta_t$ be the curve on $W_{T_0}(Y)$ obtained by following
        \begin{enumerate}
            \item the path $\gamma_t$ in $Y^h$, from $z_1=x_1(t)+ia$ to $z_2=x_2(t)+ia$;
            \item the vertical segment $\{x=x_2(t),\ a\leq y\leq a+T_0\}$ of the cylinder at $c_2$, from $w_2=x_2(t)+ia$ to $w'_2=-x_2(t)+ia$;
            \item the mirror image of $\gamma_t$ in $\overline{Y}^h$, from $z'_2=-x_2(t)+ia$ to $z'_1=-x_1(t)+ia$;
            \item the vertical segment $\{x=x_1(t),\ a\leq y\leq a+T_0\}$ of the cylinder at $c_1$, from $w'_1=-x_1(t)+ia$ to $w_1=x_1(t)+ia$.
        \end{enumerate}
        Each piece ends where the next one starts, since $w_j$ and $w'_j$ extend $z_j$ and $z'_j$ across the gluing circles, and since the mirror image of the point $z_j=x_j(t)+ia$ is the point $z'_j=-x_j(t)+ia$. Hence $\delta_t$ is a closed curve, and the curves $\delta_t$ depend continuously on $t$. For $t=0$ we have $\gamma_0=f_Y(a_e)$ and $x_j(0)=\beta_{Y,c_j}(v_j)$, so $\delta_0=f_{Y,T_0}(d_\Lambda(e))$. For $t=1$ we have $\gamma_1=e\cap Y^h$ and $x_j(1)=x_j$, so $\delta_1=d(e)$. Hence $f_{Y,T_0}(d_\Lambda(e))$ is freely homotopic to $d(e)$.
        \item The curves $d_\Lambda(e)$, $e\in\T$, are simple and pairwise disjoint, since the arcs $a_e$ are pairwise disjoint and their endpoints are pairwise distinct. We cut $\Lambda$ along all of them. Cutting $K$ along the arcs $a_e$ gives one truncated triangle $H_\tau$ for each triangle $\tau$ of $\T$, that is, a hexagon whose sides alternate between three arc sides and three intervals of $\partial K$, one for each corner of $\tau$. Cutting $\overline{K}$ gives the mirror hexagons $\overline{H}_\tau$, and cutting the annuli along the vertical segments through the endpoints of the arcs gives one rectangle $I\times[0,1]$ for each corner interval $I$, which joins $I\subset\partial H_\tau$ to its copy $\bar I\subset\partial\overline{H}_\tau$. Hence the components of the cut surface are the surfaces
        \[
            P_\tau=H_\tau\cup\overline{H}_\tau\cup(\text{three rectangles}),
        \]
        one for each triangle $\tau$ of $\T$. Write $P_\tau=A\cup B$, where $A=H_\tau\sqcup\overline{H}_\tau$ is the union of the two hexagons, $B$ is the union of the three rectangles, and $A\cap B$ is the disjoint union of the six intervals along which the rectangles are glued. Since discs and intervals have Euler characteristic $1$,
        \[
            \chi(P_\tau)=\chi(A)+\chi(B)-\chi(A\cap B)=2+3-6=-1.
        \] 
        Write $e_1,e_2,e_3$ for the sides of $\tau$ and $I_k$ for the corner interval of $H_\tau$ between $a_{e_k}$ and $a_{e_{k+1}}$, indices modulo $3$, so that the sides of $H_\tau$ are, in cyclic order, $a_{e_1},I_1,a_{e_2},I_2,a_{e_3},I_3$, see Figure~\ref{Fig-pants-tau}.      
        \begin{figure}[h]
            \centering
            \begin{overpic}[width=1.0\linewidth,keepaspectratio]{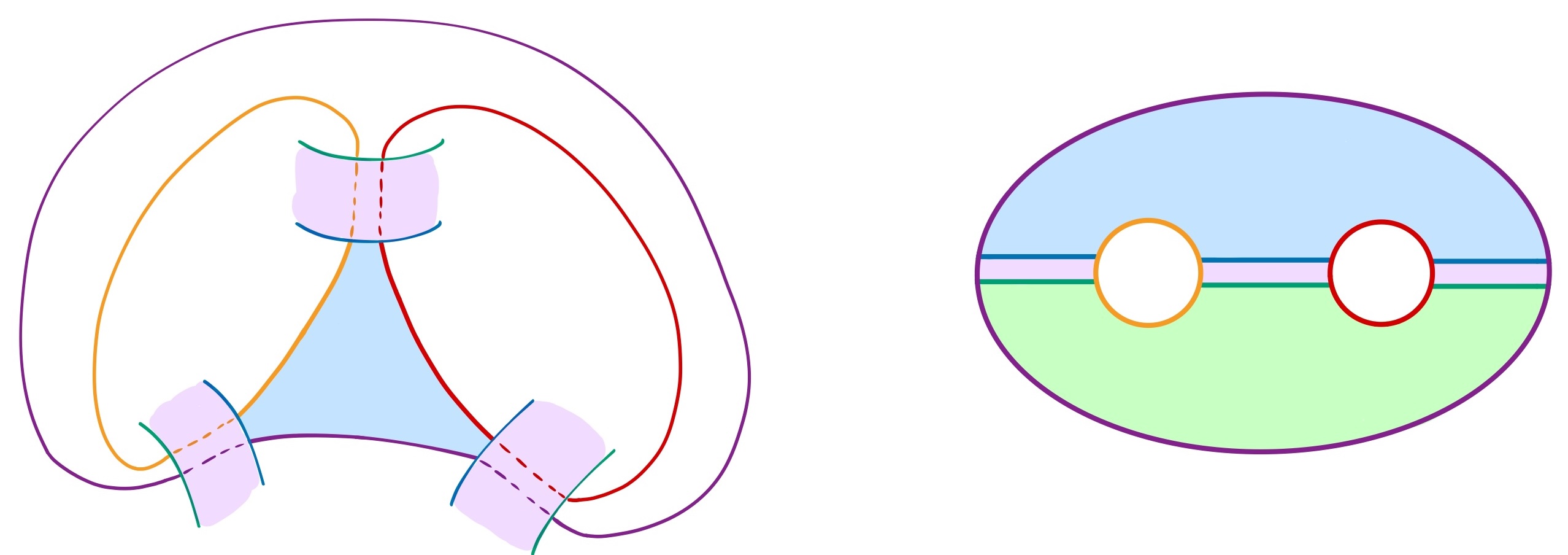}
                \put(20.5,12){\color{MidnightBlue}$H_\tau$}
                \put(28.7,19.5){\footnotesize\color{MidnightBlue}$I_2$}
                \put(28.7,25.3){\footnotesize\color{ForestGreen}$\bar I_2$}
                \put(17.3,3.3){\footnotesize\color{MidnightBlue}$I_3$}
                \put(10,0.3){\footnotesize\color{ForestGreen}$\bar I_3$}
                \put(34.5,10.3){\footnotesize\color{MidnightBlue}$I_1$}
                \put(39.5,6.5){\footnotesize\color{ForestGreen}$\bar I_1$}
                \put(42,27){\footnotesize\color{Plum}$a_{e_1}$}
                \put(27,15){\footnotesize\color{red}$a_{e_2}$}
                \put(16,16){\footnotesize\color{Orange}$a_{e_3}$}
                \put(77.5,25){\color{MidnightBlue}$H_\tau$}
                \put(77.5,8.5){\color{ForestGreen}$\overline{H}_\tau$}
                \put(64.5,20){\footnotesize\color{MidnightBlue}$I_3$}
                \put(64.5,14.5){\footnotesize\color{ForestGreen}$\bar I_3$}
                \put(80,19.6){\footnotesize\color{MidnightBlue}$I_2$}
                \put(80,14.5){\footnotesize\color{ForestGreen}$\bar I_2$}
                \put(93,19.8){\footnotesize\color{MidnightBlue}$I_1$}
                \put(93,14.5){\footnotesize\color{ForestGreen}$\bar I_1$}
                \put(76,30.8){\footnotesize\color{Plum}$a_{e_1}$}
                \put(83.5,22.3){\footnotesize\color{red}$a_{e_2}$}
                \put(68,22.3){\footnotesize\color{Orange}$a_{e_3}$}
            \end{overpic}
            \caption{The pair of pants $P_\tau$ in the proof of Lemma~\ref{lem:markings}(2), for a triangle $\tau$ with sides $e_1,e_2,e_3$.}\label{Fig-pants-tau}
        \end{figure} 
        Note that two sides of $\tau$ may lie on the same arc. In $P_\tau$, the sides $I_k$ of $H_\tau$ and $\bar I_k$ of $\overline H_\tau$ are glued to the rectangle $I_k\times[0,1]$, so the boundary of $P_\tau$ consists of the arc sides of $H_\tau$ and $\overline H_\tau$ and of the vertical sides of the rectangles. Since $a_{e_k}$ goes from an endpoint of $I_{k-1}$ to an endpoint of $I_k$, the boundary circle of $P_\tau$ through $a_{e_k}$ follows $a_{e_k}$, the vertical side of $I_k\times[0,1]$ at the common endpoint of $a_{e_k}$ and $I_k$, the mirror side $\bar a_{e_k}$ of $\overline H_\tau$, and the vertical side of $I_{k-1}\times[0,1]$ at the common endpoint of $I_{k-1}$ and $a_{e_k}$. This circle is a copy of $d_\Lambda(e_k)$. Hence $P_\tau$ has three boundary circles, and since $\chi(P_\tau)=-1$, it has genus $0$, that is, $P_\tau$ is a pair of pants.

        \begin{figure}[h]
            \centering
            \begin{overpic}[width=1.0\linewidth,keepaspectratio]{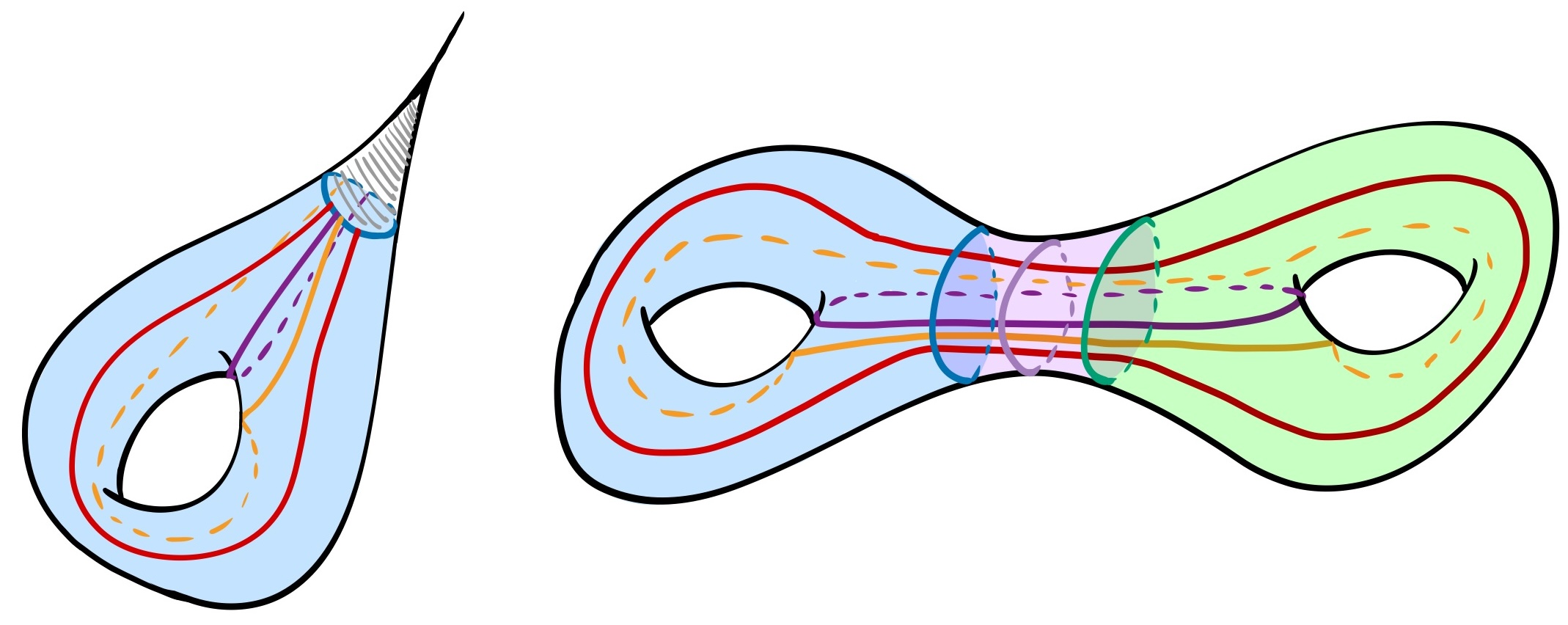}
            \end{overpic}
            \caption{A pants decomposition $\mathcal{P}_{\T}$ for a once-punctured torus $Y$ and an ideal triangulation $\T$ (arcs in red, orange, purple).}\label{Fig-pants-decomposition}
        \end{figure}

        Since every component of $\Lambda$ cut along the curves $d_\Lambda(e)$ is a pair of pants $P_\tau$, no curve $d_\Lambda(e)$ bounds a disc and no two of them bound an annulus. Hence the curves $d_\Lambda(e)$ are essential and pairwise non-homotopic, and since their complement is a union of pairs of pants, they form a pants decomposition $\mathcal P_{\T}$ of $\Lambda$. 
    \end{enumerate}
    
\end{proof}

\begin{remark}\label{rem:markings-same-arcs}
    The curves $d_\Lambda(e)$ only depend on $K$ and on the arcs $a_e$, which are fixed once and for all and used for every $Y\in\teich(\Sigma)$. Hence the pants decomposition $\mathcal P_{\T}$ of $\Lambda$ and the maximally noded surface $P_{\T}$ do not depend on $Y$.
\end{remark}

\begin{cor}\label{cor:halfdetour}\label{eq:halfdetour}
    For every $Y\in\teich(\Sigma)$,
    \begin{equation*}
        d_{\WP}^{\mathrm{aux}}\bigl(D(Y),P_{\T}\bigr)\leq2\sqrt n+\sqrt{2\pi\Bigl(2\sqrt2\,\mathcal{L}_{\T}(Y)+4E\log(1+\sqrt2)\Bigr)}.
    \end{equation*}
\end{cor}

\begin{proof}
    By Lemma~\ref{lem:markings}, the curve $f_{Y,T_0}(d_\Lambda(e))$ is freely homotopic to $d(e)$ on $W_{T_0}(Y)$. Applying~\eqref{eq:pinch} on $\teich(\Lambda)$ to $W_{T_0}(Y)$ and to the pants decomposition $\mathcal P_{\T}$ we obtain
    \begin{equation*}
         d_{\WP}^{\mathrm{aux}}\bigl(W_{T_0}(Y),P_{\T}\bigr)\leq\sqrt{2\pi\sum_{e\in\T}\ell_{W_{T_0}(Y)}(d(e))},
    \end{equation*}
    and Lemma~\ref{lem:arcs} gives
    \begin{equation}\label{eq:1}
         d_{\WP}^{\mathrm{aux}}\bigl(W_{T_0}(Y),P_{\T}\bigr)\leq\sqrt{2\pi\Bigl(2\sqrt2\,\mathcal{L}_{\T}(Y)+4E\log(1+\sqrt2)\Bigr)}.
    \end{equation}
    By Lemma~\ref{lem:opening}, for $T=T_0=2\pi$, we have 
    \begin{equation}\label{eq:2}
        d_{\WP}^{\mathrm{aux}}(W_{T_0}(Y),D(Y))\leq2\sqrt{2\pi n/T_0}=2\sqrt n.
    \end{equation} 
    Thus, by the triangle inequality, from~\eqref{eq:1} and~\eqref{eq:2} we obtain
    \begin{align*}
         d_{\WP}^{\mathrm{aux}}\bigl(D(Y),P_{\T}\bigr)&\leq d_{\WP}^{\mathrm{aux}}(D(Y),W_{T_0}(Y))+d_{\WP}^{\mathrm{aux}}\bigl(W_{T_0}(Y),P_{\T}\bigr)\\
         &\leq2\sqrt n+\sqrt{2\pi\Bigl(2\sqrt2\,\mathcal{L}_{\T}(Y)+4E\log(1+\sqrt2)\Bigr)}.
    \end{align*}
\end{proof}

\subsection{The distance in the auxiliary stratum}\label{subsec:seam-stratum}
Let $Y,Z\in\teich(\Sigma)$. The noded surfaces $D(Y)$ and $D(Z)$ lie in the stratum $S(C)\subset\overline{\teich}(\Lambda)$, where $C\subset\Lambda$ is the multicurve formed by the core curves of the annuli, whose curves are collapsed to the nodes of $D(Y)$ and $D(Z)$. We compute their distance. See also~\cite[Section~2.4]{Mondello2011}.

\begin{lemma}\label{lem:product-distance}
    For every $Y,Z\in\teich(\Sigma)$,
    \[
        d_{\WP}^{\mathrm{aux}}\bigl(D(Y),D(Z)\bigr)
        =\sqrt2\,\dwp(Y,Z).
    \]
\end{lemma}
\begin{proof}
    The surface $\Lambda\smallsetminus C$ has two components: the interior of $K$ together with the halves of the annuli adjacent to it, which is homeomorphic to $\Sigma$, and the corresponding component containing $\overline{K}$, which is homeomorphic to the mirror surface $\overline{\Sigma}$. By Section~\ref{subsec:augmentedTeich}, the stratum $S(C)$ is identified with $\teich(\Sigma)\times\teich(\overline{\Sigma})$. Since the markings $f_{Y,T}$ restrict to the markings of $Y$ and $\overline{Y}$ on $K$ and $\overline{K}$, the point $D(Y)$ corresponds to $(Y,\overline{Y})$, and similarly $D(Z)$ corresponds to $(Z,\overline{Z})$. By~\cite[Theorem~2(1), p.~17]{Wolpert2003}, the restriction of the WP metric of $\overline{\teich}(\Lambda)$ to $S(C)$ is the product of the WP metrics of $\teich(\Sigma)$ and $\teich(\overline{\Sigma})$. We denote by $d_C$ the distance on $S(C)$ for this product metric. Since $D(Y)$ and $D(Z)$ both lie in $S(C)$, their ambient distance equals their distance for the product metric (see Section~\ref{subsec:augmentedTeich}), that is
    \begin{equation}\label{eq:ambient-product}
        d_{\WP}^{\mathrm{aux}}\bigl(D(Y),D(Z)\bigr)=d_C\bigl((Y,\overline{Y}),(Z,\overline{Z})\bigr).
    \end{equation}

    Let $\sigma\colon\teich(\Sigma)\to\teich(\overline{\Sigma})$ be the map $X\mapsto\overline{X}$, where $\overline{X}$ carries the same marking as $X$. If $\zeta=x+iy$ is a holomorphic local coordinate of $X$, then $\xi=-\bar\zeta=-x+iy=u+iv$ is a holomorphic local coordinate of $\overline{X}$. Since $|d\xi|^2=|d\zeta|^2$ the surface $\overline{X}$ has the same hyperbolic metric as $X$, that is $\rho^2|d\zeta|^2=\rho^2|d\xi|^2$. Then we consider the complex conjugations 
    \begin{equation}\label{eq:conj-maps}
        \begin{aligned}
            Q(X)&\longrightarrow Q(\overline{X}), &\qquad \varphi=\varphi(\zeta)\,d\zeta^2&\longmapsto\overline{\varphi}=\overline{\varphi(\zeta)}\,d\bar\zeta^2=\overline{\varphi(-\bar\xi)}\,d\xi^2,\\
            B(X)&\longrightarrow B(\overline{X}), &\qquad \mu=\mu(\zeta)\,\frac{d\bar\zeta}{d\zeta}&\longmapsto\overline{\mu}=\overline{\mu(\zeta)}\,\frac{d\zeta}{d\bar\zeta}=\overline{\mu(-\bar\xi)}\,\frac{d\bar\xi}{d\xi}.
        \end{aligned}
    \end{equation}

    The first map is well-defined since $\xi\mapsto\overline{\varphi(-\bar\xi)}$ is holomorphic, and it is a bijection, whose inverse is the same construction from $\overline{X}$ to $\overline{\overline{X}}=X$.
    Then, since $d\xi=-d\bar\zeta$, in the coordinate $\xi$ of $\overline{X}$ we have $ \overline{\varphi}=\overline{\varphi(\zeta)}\,d\xi^2$ and $\overline{\mu}=\overline{\mu(\zeta)}\,\frac{d\bar\xi}{d\xi}$. 
    Therefore
    \begin{equation}\label{eq:conj-pairing}
            \int_{\overline{X}}\overline{\mu}\,\overline{\varphi}
            =\int_{\overline{X}}\overline{\mu(-\bar\xi)}\;\overline{\varphi(-\bar\xi)}\,du\,dv
            =\int_{X}\overline{\mu(\zeta)\,\varphi(\zeta)}\,dx\,dy
            =\overline{\int_{X}\mu\varphi},
    \end{equation}
    and 

    \begin{equation}\label{eq:conj-norm}
        \begin{aligned}
            \|\overline{\varphi}\|_2^2
            &=\int_{\overline{X}}\bigl|\overline{\varphi(-\bar\xi)}\bigr|^2\rho(-\bar\xi)^{-2}\,du\,dv
            =\int_{\overline{X}}|\varphi(-\bar\xi)|^2\rho(-\bar\xi)^{-2}\,du\,dv\\
            &=\int_{X}|\varphi(\zeta)|^2\rho(\zeta)^{-2}\,dx\,dy=\|\varphi\|_2^2.
        \end{aligned}
    \end{equation}
    The differential of $\sigma$ maps $[\mu]$ to $[\overline{\mu}]$. Indeed, let $\mu\in B(X)$ and let $(X_t)_t$ be a smooth path in $\teich(\Sigma)$ with $X_0=X$, described by Beltrami coefficients $\mu(t)$ on $X$ with $\mu(0)=0$ and $\dot\mu(0)=\mu$, so that the tangent vector of the path at $t=0$ is $[\mu]$. If $w$ is a local holomorphic coordinate of $X_t$, viewed as a function of $\zeta$, then $\partial_{\bar\zeta}w=\mu(t)\,\partial_\zeta w$. The function $\hat w=-\bar w$ is a local holomorphic coordinate of $\overline{X_t}$. Since $\zeta=-\bar\xi$, we have $\partial_{\bar\xi}=-\partial_\zeta$ and $\partial_\xi=-\partial_{\bar\zeta}$, and for every function $f$ we have $\partial_\zeta\bar f=\overline{\partial_{\bar\zeta}f}$ and $\partial_{\bar\zeta}\bar f=\overline{\partial_\zeta f}$. Hence
    \[
        \partial_{\bar\xi}\hat w=-\partial_\zeta(-\bar w)=\overline{\partial_{\bar\zeta}w}=\overline{\mu(t)}\;\overline{\partial_\zeta w}=\overline{\mu(t)}\;\partial_{\bar\zeta}\bar w=\overline{\mu(t)}\;\partial_\xi\hat w.
    \]
    Thus $\sigma(X_t)=\overline{X_t}$ is described by the Beltrami coefficients $\overline{\mu(t)}$ on $\overline{X}$, whose derivative at $t=0$ is $\overline{\mu}$. Hence the tangent vector of the path $(\sigma(X_t))_t$ at $t=0$ is $[\overline{\mu}]$, that is, $d\sigma([\mu])=[\overline{\mu}]$. This is well defined, since $\mu\in N(X)$ if and only if $\overline{\mu}\in N(\overline{X})$ by~\eqref{eq:conj-pairing}.

    Consequently the map $\sigma$ preserves the WP norm: for every $\mu\in B(X)$,
    \[ 
        \begin{aligned}
            \bigl\|d\sigma([\mu])\bigr\|_{\WP}=\bigl\|[\overline{\mu}]\bigr\|_{\WP}
            &=\sup_{\psi\in Q(\overline{X})\smallsetminus\{0\}}\frac{\bigl|\int_{\overline{X}}\overline{\mu}\,\psi\bigr|}{\|\psi\|_2}
            \underset{\eqref{eq:conj-maps}}{=}\sup_{\varphi\in Q(X)\smallsetminus\{0\}}\frac{\bigl|\int_{\overline{X}}\overline{\mu}\,\overline{\varphi}\bigr|}{\|\overline{\varphi}\|_2}\\
            &\underset{\eqref{eq:conj-pairing},\,\eqref{eq:conj-norm}}{=}\sup_{\varphi\in Q(X)\smallsetminus\{0\}}\frac{\bigl|\int_X\mu\varphi\bigr|}{\|\varphi\|_2}
            =\bigl\|[\mu]\bigr\|_{\WP}.
        \end{aligned} 
    \]
    hence, since $\sigma$ is a bijection with inverse $\overline{X}\mapsto X$, it is an isometry. In particular $\dwp(\overline{Y},\overline{Z})=\dwp(Y,Z)$. Set $d=\dwp(Y,Z)$.

    Now we compute $d_C\bigl((Y,\overline{Y}),(Z,\overline{Z})\bigr)$. Let $\gamma(t)=(\alpha(t),\beta(t))$, $t\in[0,1]$, be a piecewise smooth path in $S(C)$ from $(Y,\overline{Y})$ to $(Z,\overline{Z})$. The product metric on $S(C)=\teich(\Sigma)\times\teich(\overline{\Sigma})$ gives to the tangent vector $(\dot\alpha,\dot\beta)$ the norm $\bigl(\|\dot\alpha\|_{\WP}^2+\|\dot\beta\|_{\WP}^2\bigr)^{1/2}$. The path $\alpha$ joins $Y$ to $Z$ in $\teich(\Sigma)$, so $\ell(\alpha)\geq d$. The path $\beta$ joins $\overline{Y}$ to $\overline{Z}$ in $\teich(\overline{\Sigma})$, so $\sigma^{-1}\circ\beta$ joins $Y$ to $Z$, and since $\sigma$ is an isometry, $\ell(\beta)=\ell(\sigma^{-1}\circ\beta)\geq d$. Write $A=\|\dot\alpha\|_{\WP}$ and $B=\|\dot\beta\|_{\WP}$. Since $2(A^2+B^2)-(A+B)^2=(A-B)^2\geq0$, we have $\sqrt{A^2+B^2}\geq(A+B)/\sqrt2$, hence
    \[
        \ell(\gamma)=\int_0^1\sqrt{A^2+B^2}\,dt\geq\frac{1}{\sqrt2}\left(\int_0^1A\,dt+\int_0^1B\,dt\right)=\frac{\ell(\alpha)+\ell(\beta)}{\sqrt2}\geq\frac{2d}{\sqrt2}=\sqrt2\,d.
    \]
    Taking the infimum over $\gamma$ gives
    \begin{equation}\label{eq:lower-bound}
        d_C\bigl((Y,\overline{Y}),(Z,\overline{Z})\bigr)\geq\sqrt2\,d.
    \end{equation}
    Let $\eta>0$, and let $\alpha$ be a piecewise smooth path from $Y$ to $Z$ in $\teich(\Sigma)$ with $\ell(\alpha)<d+\eta$, which exists by definition of $d$. The path $\gamma(t)=\bigl(\alpha(t),\sigma(\alpha(t))\bigr)$ joins $(Y,\overline{Y})$ to $(Z,\overline{Z})$ in $S(C)$. Its tangent vector is $\bigl(\dot\alpha,d\sigma(\dot\alpha)\bigr)$, and $\|d\sigma(\dot\alpha)\|_{\WP}=\|\dot\alpha\|_{\WP}$, so its norm is
    \[ 
        \|\dot\gamma\|=\sqrt{\|\dot\alpha\|_{\WP}^2+\|d\sigma(\dot\alpha)\|_{\WP}^2}=\sqrt2\,\|\dot\alpha\|_{\WP}\]
    and
    \[
        \ell(\gamma)=\sqrt2\,\ell(\alpha)<\sqrt2\,(d+\eta).
    \]
    Hence $d_C\bigl((Y,\overline{Y}),(Z,\overline{Z})\bigr)<\sqrt2\,(d+\eta)$ for every $\eta>0$, that is,
    \begin{equation}\label{eq:upper-bound}
        d_C\bigl((Y,\overline{Y}),(Z,\overline{Z})\bigr)\leq\sqrt2\,d.
    \end{equation}
    By~\eqref{eq:lower-bound} and~\eqref{eq:upper-bound}, we obtain $d_C\bigl((Y,\overline{Y}),(Z,\overline{Z})\bigr)=\sqrt2\,d$, and by~\eqref{eq:ambient-product} this gives $d_{\WP}^{\mathrm{aux}}\bigl(D(Y),D(Z)\bigr)=\sqrt2\,\dwp(Y,Z)$, which concludes the proof.
\end{proof}

\subsection{Proof of the main theorem}\label{subsec:proof-main}

\begin{theorem}\label{mainthm}
Let $g,n$ be integers with $n\geq1$ and $2g-2+n>0$. There is a universal constant $C<6$ such that any two hyperbolic surfaces $Y,Z\in\teich_{g,n}$, equipped with the same ideal triangulation $\T$, satisfy
\[
    \dwp(Y,Z)
    \le C\sqrt{n+\mathcal{L}_{\T}(Y)+\mathcal{L}_{\T}(Z)}.
\]
\end{theorem}

\begin{proof}
    Let $Y,Z\in\teich_{g,n}$ be equipped with the same topological ideal triangulation $\T$, and mark $W_{T_0}(Y)$ and $W_{T_0}(Z)$ by $f_{Y,T_0}$ and $f_{Z,T_0}$, see Section~\ref{subsec:opening}.

    By Lemma~\ref{lem:product-distance},
    \begin{equation}\label{eq:main-step1}
        \sqrt2\,\dwp(Y,Z)=d_{\WP}^{\mathrm{aux}}\bigl(D(Y),D(Z)\bigr).
    \end{equation}

    By Lemma~\ref{lem:markings} and Remark~\ref{rem:markings-same-arcs}, the maximally noded surface $P_{\T}\in\overline{\teich}(\Lambda)$ is the same for $Y$ and $Z$. By the triangle inequality,
    \begin{equation}\label{eq:main-step2}
        d_{\WP}^{\mathrm{aux}}\bigl(D(Y),D(Z)\bigr)\leq d_{\WP}^{\mathrm{aux}}\bigl(D(Y),P_{\T}\bigr)+d_{\WP}^{\mathrm{aux}}\bigl(P_{\T},D(Z)\bigr).
    \end{equation}

   Let $V\in\{Y,Z\}$. By~\eqref{eq:delta-lower}, every arc $e$ of $\T$ satisfies $\delta_e(V)\geq2\log4$, hence $\mathcal{L}_{\T}(V)\geq2E\log4$, that is, $E\leq\mathcal{L}_{\T}(V)/(2\log4)$. Therefore
    \[
        2\sqrt2\,\mathcal{L}_{\T}(V)+4E\log(1+\sqrt2)\leq\Bigl(2\sqrt2+\frac{2\log(1+\sqrt2)}{\log4}\Bigr)\mathcal{L}_{\T}(V)=B\,\mathcal{L}_{\T}(V),
    \]
    where $B=2\sqrt2+\log(1+\sqrt2)/\log2$, and Corollary~\ref{cor:halfdetour} gives
    \begin{equation}\label{eq:main-step3}
        d_{\WP}^{\mathrm{aux}}\bigl(D(V),P_{\T}\bigr)\leq2\sqrt n+\sqrt{2\pi B\,\mathcal{L}_{\T}(V)}.
    \end{equation}

    Combining~\eqref{eq:main-step1}, \eqref{eq:main-step2} and~\eqref{eq:main-step3} for $V=Y$ and $V=Z$, and dividing by $\sqrt2$, we obtain
    \[
        \begin{aligned}
            \dwp(Y,Z)&\leq\frac{1}{\sqrt2}\Bigl(4\sqrt n+\sqrt{2\pi B}\bigl(\sqrt{\mathcal{L}_{\T}(Y)}+\sqrt{\mathcal{L}_{\T}(Z)}\bigr)\Bigr)\\
            &=2\sqrt2\,\sqrt n+\sqrt{\pi B}\Bigl(\sqrt{\mathcal{L}_{\T}(Y)}+\sqrt{\mathcal{L}_{\T}(Z)}\Bigr).
        \end{aligned}
    \]
    By the Cauchy--Schwarz inequality,
    \[
        \dwp(Y,Z)\leq\sqrt{8+2\pi B}\,\sqrt{n+\mathcal{L}_{\T}(Y)+\mathcal{L}_{\T}(Z)}.
    \]
    Finally $B<4.1$, so $C=\sqrt{8+2\pi B}<5.811<6$. This constant does not depend on $g$ and $n$, on the number of arcs at each cusp, nor on the geometry of $Y$ and $Z$.
\end{proof}

We now deduce Corollary~\ref{mainthm2} from Theorem~\ref{mainthm} and Lemma~\ref{lem:shear-trunc}.

\begin{cor}\label{mainthm2}
Let $g,n$ be integers with $n\geq1$ and $2g-2+n>0$. Any two hyperbolic surfaces $Y,Z\in\teich_{g,n}$, equipped with the same topological ideal triangulation $\T$, satisfy
\[
    \dwp(Y,Z)\leq 21\sqrt{(2g-2+n)\Bigl(\log\bigl(24(2g-2+n)\bigr)+\mathcal{O}_{\T}(Y)+\mathcal{O}_{\T}(Z)\Bigr)},
\]
where $\mathcal{O}_{\T}$ denotes the oscillation norm of Kahn--Markovi\'c~\cite{KahnMarkovic2008}, defined in Section~\ref{subsec:shear}.
\end{cor}

\begin{proof}[Proof of Corollary~\ref{mainthm2}]
    Set $m=2g-2+n$, so that $m\geq1$ and $E=6g-6+3n=3m$. For $V\in\{Y,Z\}$, Lemma~\ref{lem:shear-trunc} gives
    \[
        \mathcal{L}_{\T}(V)\leq2E\bigl(\log(8E)+\mathcal{O}_{\T}(V)\bigr)=6m\bigl(\log(24m)+\mathcal{O}_{\T}(V)\bigr).
    \]
    Moreover $n=m+2-2g\leq m+2\leq3m$, since $m\geq1$, and $3<\log24\leq\log(24m)$, so $n\leq m\log(24m)$. Therefore
    \[
        \begin{aligned}
            n+\mathcal{L}_{\T}(Y)+\mathcal{L}_{\T}(Z)&\leq13\,m\log(24m)+6m\bigl(\mathcal{O}_{\T}(Y)+\mathcal{O}_{\T}(Z)\bigr)\\
            &\leq13\,m\bigl(\log(24m)+\mathcal{O}_{\T}(Y)+\mathcal{O}_{\T}(Z)\bigr).
        \end{aligned}
    \]
    By Theorem~\ref{mainthm}, with $C=\sqrt{8+2\pi B}<5.811$,
    \[
        \dwp(Y,Z)\leq C\sqrt{13}\,\sqrt{m\bigl(\log(24m)+\mathcal{O}_{\T}(Y)+\mathcal{O}_{\T}(Z)\bigr)},
    \]
    and $C\sqrt{13}<5.811\cdot3.606<21$.
\end{proof}

\begin{cor}\label{cor:shear-l1}
Let $g,n$ be integers with $n\geq1$ and $2g-2+n>0$.
Any two surfaces $Y,Z\in\teich_{g,n}$, equipped with the
same marked ideal triangulation $\T$, satisfy
\[
    \dwp(Y,Z)
    \leq 21\sqrt{(2g-2+n)\left(\log(24(2g-2+n))
             +\sum_{e\in\T}\bigl(|s_e(Y)|+|s_e(Z)|\bigr)\right)}.
\]
\end{cor}

\begin{proof}
Fix $X\in\teich_{g,n}$ and a cusp $c$. Write
$a_r=s_{e_{c,r}}(X)$ for $1\leq r\leq v(c)$, using the
endpoint notation of Section~\ref{subsec:shear}.
Since the arcs are ideal one has $\sum_{r=1}^{v(c)}a_r=0$. Hence, for
every consecutive set $I$ of endpoints at $c$,
\[
    \left|\sum_{r\in I}a_r\right|
    =\frac12\left|\sum_{r\in I}a_r-\sum_{r\notin I}a_r\right|
    \leq\frac12\sum_{r=1}^{v(c)}|a_r|
    \leq\sum_{e\in\T}|s_e(X)|.
\]
The last inequality holds because each arc has at most two
endpoints at $c$. Taking the maximum over $c$ and $I$ gives
$\mathcal{O}_{\T}(X)\leq\sum_{e\in\T}|s_e(X)|$ and we applly this inequality to $Y$ and $Z$ in Corollary~\ref{mainthm2} to conclude.
\end{proof}

\bibliographystyle{alpha}
\bibliography{bib}

\end{document}

%% file: commands.tex
\newtheorem{theorem}{Theorem}[section]
\newtheorem{proposition}[theorem]{Proposition}
\newtheorem{lemma}[theorem]{Lemma}
\newtheorem*{lemma*}{Lemma}
\newtheorem*{claim*}{Claim}
\newtheorem*{definition*}{Definition}
\newtheorem{cor}[theorem]{Corollary}
\theoremstyle{definition}

\newtheorem{claim}[theorem]{Claim}

\newtheorem{remark}[theorem]{Remark}

\newtheorem*{theorem**}{Theorem\theoremnum}

\newenvironment{theorem*}[1][]{%
  \if\relax\detokenize{#1}\relax
    \begin{theorem}
  \else
    \renewcommand{\thetheorem}{#1}
    \begin{theorem}%
    \addtocounter{theorem}{-1}
  \fi
}{%
  \end{theorem}
}

\newenvironment{cor*}[1][]{%
  \if\relax\detokenize{#1}\relax
    \begin{cor}
  \else
    \renewcommand{\thetheorem}{#1}
    \begin{cor}%
    \addtocounter{theorem}{-1}
  \fi
}{%
  \end{cor}
}

\newtheorem*{proposition**}{Proposition\theoremnum}

\newenvironment{proposition*}[1][]{%
  \if\relax\detokenize{#1}\relax
    \begin{proposition}%
  \else
    \renewcommand{\thetheorem}{#1}
    \begin{proposition}%
    \addtocounter{theorem}{-1}
  \fi
}{%
  \end{proposition}
}

\newcommand{\addQEDstyle}[2]{\AtBeginEnvironment{#1}{\pushQED{\qed}\renewcommand{\qedsymbol}{#2}}\AtEndEnvironment{#1}{\popQED}}

\addQEDstyle{example}{$\triangle$}
\newcommand{\T}{\mathcal T}

\DeclareMathOperator{\ima}{Im}

\DeclareMathOperator{\h}{\mathfrak{h}}

\DeclareMathOperator{\Mod}{Mod}

\DeclareMathOperator{\area}{area}

\DeclareMathOperator{\dwp}{d_{WP}}

\DeclareMathOperator{\WP}{WP}

\DeclareMathOperator{\teich}{Teich}

\DeclareMathOperator{\FN}{FN}

\newcommand{\C}{\mathbb{C}}

\newcommand{\N}{\mathbb{N}}
\newcommand{\R}{\mathbb{R}}
\newcommand{\Hh}{\mathbb{H}}
\newcommand{\Z}{\mathbb{Z}}

\newcommand{\Ms}{\mathcal{M}}

\newcommand{\Ps}{\mathscr{P}}

\newcommand{\Zero}{\mathbf{0}}

\newcommand{\coloneq}{\mathrel{\resizebox{\widthof{$\mathord{=}$}}{\height}{ $\!\!\resizebox{1.2\width}{0.8\height}{\raisebox{0.23ex}{$\mathop{:}$}}\!\!=\!\!$ }}}